\documentclass[a4paper]{article}
\usepackage{amsmath}
\usepackage{amssymb}
\usepackage{amsthm}
\usepackage[dvipdfmx]{graphicx}
\usepackage[dvipdfmx]{pict2e}
\usepackage{bm} 
\usepackage{color}
\newtheorem{theorem}{Theorem}[section]
\newtheorem{proposition}[theorem]{Proposition}

\newtheorem{lemma}[theorem]{Lemma}
\theoremstyle{definition}
\newtheorem{definition}[theorem]{Definition}

\newtheorem{remark}[theorem]{Remark}

\newcommand{\oline}[1]{\mathbin{\overline{#1}}}
\newcommand{\uline}[1]{\mathbin{\underline{#1}}}
\begin{document}
\title{Cocycles of $G$-Alexander biquandles and $G$-Alexander multiple conjugation biquandles}
\author{Atsushi Ishii, Masahide Iwakiri, Seiichi Kamada，\\
Jieon Kim, Shosaku Matsuzaki, and Kanako Oshiro}
\date{\today}
\maketitle
\begin{abstract} 
Biquandles and multiple conjugation biquandles are algebras which are related to links and handlebody-links in $3$-space. 
Cocycles of them can be used to construct state-sum type invariants of links and handlebody-links. In this paper we discuss cocycles of a certain class of biquandles and multiple conjugation biquandles, which we call 
$G$-Alexander biquandles and $G$-Alexander multiple conjugation biquandles,  with a relationship with group cocycles.  We give a  
method to obtain a (biquandle or multiple conjugation biquandle) cocycle of them from a group cocycle. 
\end{abstract}

\section{Introduction}

A quandle  is an algebra related to knots and links (cf. \cite{FennRourke92, Joyce82, Matveev82}). Using 
a (quandle) 2-cocycle or 3-cocycle in the (co)homology theory of a quandle (cf. \cite{CarterJelsovskyKamadaLangfordSaito03, FennRourkeSanderson95, FennRourkeSanderson2007}), we can construct an invariant of knots and links.
A multiple conjugation quandle is introduced in \cite{Ishii15MCQ}. 
It is an algebra motivated by study of handlebody-links, which are isotopy classes of embedded handlebodies in 3-space. 
Another notion called a qualgebra was introduced in \cite{Lebed14, Lebed15} which is related to trivalent graphs in 3-space. 

The notions of a quandle and a multiple conjugation quandle are generalized to biquandles (cf. \cite{JSC_ME_MS2004, CenicerosElhamdadiGreenNelson14, FennJordanKauffman04, KauffmanRadford03}) and multiple conjugation biquandles (cf. \cite{IshiiIwakiriKamadaKimMatsuzakiOshiroPMB, IshiiNelson16}), although the terminologies and definitions of such generalizations vary in the literature.    
Homology theory on biquandles is found in \cite{JSC_ME_MS2004, FennRourkeSanderson2007} and homology theory on multiple conjugation biquandles is found in \cite{IIKKMO2}. In \cite{IIKKMO2} we discussed 
the (co)homology theory of multiple conjugation biquandles and we provided a method of constructing an invariant of handlebody-links from a 2 or 3-cocycle. Thus, in order to construct such an invariant of links or handlebody-links, it is important to find a cocycle of a biquandle or a multiple conjugation biquandle.  
In this paper, we focus on a certain class of biquandles and multiple conjugation biquandles which we call $G$-Alexander biquandles and $G$-Alexander multiple conjugation biquandles.  

The quandles associated to  $G$-families of Alexander quandles \cite{IshiiIwakiriJangOshiro13} are an important class of quandles. In \cite{Nosaka13}, Nosaka studied  cocycles of a $G$-family $M$ of Alexander quandles, for a group $G$ and a $G$-module $M$, and gave a method of obtaining a cocycle of the associated quandle  $M \times G$. 
For a given $G$-invariant group cocycle of the $G$-family $M$, we can construct a cocycle of the quandle $M \times G$ by using invariant theory. Our approach in this paper is an analogue of Nosaka's argument. 

The paper is organized as follows. In Section \ref{Preliminaries}, we review 
definitions related to the (co)homology of biquandles, multiple conjugation biquandles and groups. 
We also give the definition of $G$-Alexander biquandles and $G$-Alexander multiple conjugation biquandles.  
In Sections \ref{app:11} and \ref{app:12}
we focus on $G$-Alexander biquandles $X = M\times G$ and cocycles of them. The case where the biquandle $X$ has the trivial $X$-set (or has $X$ itself as an $X$-set, resp.) is discussed in Section~\ref{app:11} (or \ref{app:12}).   
We show how to obtain a (birack) cocycle of the biquandle $X= M\times G$ from a group cocycle: Theorems~\ref{birack_cocycle_main} and~\ref{birack cocycle explicit}  (Theorems~\ref{birack_cocycle_with_X-set_main} and~\ref{birack_cocycle_explicit_Xset}).   
In Sections \ref{app:21} and \ref{app:22}
we focus on $G$-Alexander multiple conjugation biquandles $X =  \bigsqcup_{m \in M}(\{m\} \times G) = M\times G$ and cocycles of them. The case where the biquandle $X$ has the trivial $X$-set (or has $X$ itself as an $X$-set, resp.) is discussed in Section~\ref{app:21} (or \ref{app:22}).   
We show how to obtain a (multiple conjugation biquandle) cocycle of the $G$-Alexander multiple conjugation biquandle from a group cocycle: Theorems~\ref{mcb_cocycles_no_X-set_main} and~\ref{mcb_cocycle_no_X-setexplicit}  
(Theorems~\ref{mcb_cocycle_with_X-set_main},    \ref{mcb_cocycle_with_X-set_main_explicit1} and \ref{mcb_cocycle_with_X-set_main_explicit2}).  
Since any $G$-invariant multilinear map provides a $G$-invariant group cocycle, we can obtain many cocycles.

This research is supported by JSPS KAKENHI Grant Numbers  
16K17600, 
18K03292 
and
19H01788.  
It is also supported by 
Young Researchers Program through the National Research Foundation of Korea
(NRF) funded by the Ministry of Education, Science and Technology (NRF-2018R1C1B6007021).

\section{Preliminaries}\label{Preliminaries}

In this section, we review definitions of biquandles and multiple conjugation biquandles, and the (co)homology of 
those and groups. 
In Subsection~\ref{subsect_G_Alexander}, we give the definition of $G$-Alexander biquandles and $G$-Alexander multiple conjugation biquandles.  

\subsection{Biquandles and multiple conjugation biquandles (MCB)} 
\begin{definition}[\cite{FennRourkeSanderson95,KauffmanRadford03}]
A \textit{biquandle} is a non-empty set $X$ with binary operations $\uline{*},\oline{*}:X\times X\to X$ satisfying the following axioms.
\begin{itemize}
\item[(B1)]
For any $x\in X$, $x\uline{*}x=x\oline{*}x$.
\item[(B2)]
For any $a\in X$, the map $\uline{*}a:X\to X$ sending $x$ to $x\uline{*}a$ is bijective.
\item[]
For any $a\in X$, the map $\oline{*}a:X\to X$ sending $x$ to $x\oline{*}a$ is bijective.
\item[]
The map $S:X\times X\to X\times X$ defined by $S(x,y)=(y\oline{*}x,x\uline{*}y)$ is bijective.
\item[(B3)]
For any $x,y,z\in X$,
\begin{align*}
&(x\uline{*}y)\uline{*}(z\uline{*}y)=(x\uline{*}z)\uline{*}(y\oline{*}z), \\
&(x\uline{*}y)\oline{*}(z\uline{*}y)=(x\oline{*}z)\uline{*}(y\oline{*}z), \\
&(x\oline{*}y)\oline{*}(z\oline{*}y)=(x\oline{*}z)\oline{*}(y\uline{*}z).
\end{align*}
\end{itemize}
\end{definition}

A \textit{birack} is a non-empty set $X$ with binary operations $\uline{*},\oline{*}:X\times X\to X$ satisfying (B2) and (B3). 

\begin{definition}
For a biquandle (or a birack) $X$, an \textit{$X$-set} is a non-empty set $Y$ with a map $*:Y\times X\to Y$ satisfying the following axioms.
\begin{itemize}
\item
For any $y\in Y$ and $a,b\in X$, $(y*a)*(b\oline{*}a)=(y*b)*(a\uline{*}b)$.
\item
For any $x \in X$, the map $*x:Y \to Y; y \mapsto y*x $ is bijective.
\end{itemize}
\end{definition}


Let $X$ be the disjoint union of groups $G_\lambda$ ($\lambda\in\Lambda$). We denote by $G_a$ the group $G_\lambda$ to which $a\in X$ belongs. We denote by $e_\lambda$ the identity of $G_\lambda$.

\begin{definition}[\cite{IshiiIwakiriKamadaKimMatsuzakiOshiroPMB}] \label{def:MGQ2}
A {\it multiple conjugation biquandle} is a non-empty set $X$ which is 
the disjoint union of groups $G_\lambda$ ($\lambda\in\Lambda$) with binary operations $\uline{*},\oline{*}:X\times X\to X$ 
satisfying the following axioms. 
\begin{itemize}
\item
For any $x,y,z\in X$,
\begin{align}
&(x\uline{*}y)\uline{*}(z\uline{*}y)=(x\uline{*}z)\uline{*}(y\oline{*}z), \label{eq:R3-1'} \\
&(x\uline{*}y)\oline{*}(z\uline{*}y)=(x\oline{*}z)\uline{*}(y\oline{*}z), \label{eq:R3-2'} \\
&(x\oline{*}y)\oline{*}(z\oline{*}y)=(x\oline{*}z)\oline{*}(y\uline{*}z). \label{eq:R3-3'}
\end{align}
\item
For any $a,x\in X$, the maps $\uline{*}a:G_x\to G_{x\uline{*}a}$ and $\oline{*}a:G_x\to G_{x\oline{*}a}$ are group homomorphisms.
\item
For any $a,b\in G_\lambda$ and $x\in X$,
\begin{align}
&x\uline{*}ab=(x\uline{*}a)\uline{*}(b\oline{*}a),
\hspace{1em}x\uline{*}e_\lambda=x, \label{eq:x*u(ab)'} \\
&x\oline{*}ab=(x\oline{*}a)\oline{*}(b\oline{*}a),
\hspace{1em}x\oline{*}e_\lambda=x, \label{eq:x*o(ab)'} \\
&a^{-1}b\oline{*}a=ba^{-1}\uline{*}a. \label{eq:R14'}
\end{align}
\end{itemize}
\end{definition}

A multiple conjugation biquandle $(X = \bigsqcup_{\lambda\in\Lambda}G_\lambda,  \uline*, \oline*)$ is regarded as a biquandle 
$(X,  \uline*, \oline*)$ by forgetting the decomposition $X = \bigsqcup_{\lambda\in\Lambda}G_\lambda$ and the group structure of each $G_\lambda$. 

\begin{definition}
For a multiple conjugation biquandle $X=\bigsqcup_{\lambda\in\Lambda}G_\lambda$, an \textit{$X$-set} is a non-empty set $Y$ with a map $*:Y\times X\to Y$ satisfying the following axioms.
\begin{itemize}
\item
For any $y\in Y$ and $a,b \in G_\lambda$,
$y*e_\lambda=y$, $y*(ab)=(y*a)*(b\oline{*}a)$, where $e_\lambda$ is the identity of $G_\lambda$.
\item
For any $y\in Y$ and $a,b\in X$, $(y*a)*(b\oline{*}a)=(y*b)*(a\uline{*}b)$.
\end{itemize}
\end{definition}

Any singleton set $\{y_0\}$ is also an $X$-set with the map $*$ defined by $y_0*x=y_0$ for $x\in X$, which is called a \textit{trivial $X$-set}. Any multiple conjugation biquandle $X$ itself is an $X$-set with the map $\uline{*}$ or $\oline{*}$.

\subsection{$G$-Alexander biquandles and $G$-Alexander multiple conjugation biquandles}\label{subsect_G_Alexander}

\begin{definition}\label{def_Gfamily_of_Alexander_quandle}
Let $R$ be a ring, $G$ a group and $M$ a right $R[G]$-module, where $R[G]$ is the group ring of $G$ over $R$. Let $\varphi: G \to Z(G)$ be a group homomorphism, where $Z(G)$ denotes the center of $G$. We define binary operations $\uline{*}$ and $\oline{*} : (M \times G)^2\to M \times G$ by
\begin{align*}
(m,g)\uline{*}(n,h) &:=(mh+n(\varphi(h)-h), h^{-1}gh),\\
(m,g)\oline{*}(n,h) &:=(m\varphi(h), g).
\end{align*}
Then, $(M \times G, \uline*, \oline*)$ is a biquandle (cf.\cite{IshiiNelson16}), which we call the {\it $G$-Alexander biquandle} of $(M, \varphi)$, and $(M \times G=\bigsqcup_{m \in M}(\{m\} \times G), \uline*, \oline*)$ is a multiple conjugation biquandle by 
$$(m,g) (m,h) := (m, gh)$$
for any $m \in M$ and any $g, h \in G$ (cf.\cite{IshiiNelson16}), which we call the {\it $G$-Alexander multiple conjugation biquandle} of $(M, \varphi)$. 
\end{definition}

Unless otherwise stated, we suppose that the ring $R$ is $\mathbb Z$. 

\begin{remark}
In general, for a given $G$-family of biquandles, one can associate a biquandle   
called the {\it associated biquandle of the $G$-family of biquandles} (\cite{IshiiNelson16}), and 
  a multiple conjugation biquandle  
called the {\it associated multiple conjugation biquandle of the $G$-family of biquandles} or 
the {\it partially multiplicative biquandle associated to the $G$-family of biquandles} 
(\cite{IshiiNelson16}).  
The $G$-Alexander biquandle (or the $G$-Alexander multiple conjugation biquandle) in the definition above is the associated biquandle (or the associated multiple conjugation biquandle) of a \lq\lq  $G$-family of Alexander biquandle\rq\rq \,  in the sense of \cite{IshiiNelson16}. 
\end{remark} 

\subsection{The (co)homology of biquandles}\label{subsec_homology_biquandles}
In this subsection, we review the birack chain complex for a biquandle.

Let $X$ be a biquandle (or more generally a birack) and $Y$ an $X$-set. Let $C^{\rm BR}_n( X ; \mathbb Z)_Y$  denote the free abelian group generated by the elements $\bm{x}=(y, x_1, \ldots, x_n) \in Y \times X^n$ if $n\geq 1$, and $C^{\rm BR}_n(X ; \mathbb Z)_Y=0$ otherwise. 
(The superscripts BR means \textit{birack}.)

For integers $n$ and $i$ such that $n \ge 2$ and $1 \le i \le n$ and for any element $\bm{x} =(y, x_1, \ldots, x_i, \ldots,  x_n) \in Y \times X^n$, we write
\begin{align*}
(\bm{x}^{i-1}, \bm{x}_{i+1}) &:= (y, x_1, \ldots, x_{i-1},x_{i+1}, \ldots, x_n) \text{ and}\\
(\bm{x}^{i-1} \uline* x_i, \bm{x}_{i+1} \oline* x_i) &:=(y*x, x_1 \uline* x_i, \ldots, x_{i-1} \uline* x_i,x_{i+1} \oline* x_i, \ldots, x_n \oline* x_i) ,
\end{align*}
where these elements belong to $Y \times X^{n-1}$.

Define a boundary map $\partial_n^{\rm BR} : C^{\rm BR}_n(X ; \mathbb Z)_Y \to C^{\rm BR}_{n-1} (X ; \mathbb Z)_Y $ by 
\[
\partial_n^{\rm BR} (\bm{x} ) := \displaystyle \sum_{i=1}^{n} (-1)^{i-1} \Bigl\{ (\bm{x}^{i-1} ,\bm{x}_{i+1})-( \bm{x}^{i-1} \uline* x_i ,\bm{x}_{i+1} \oline* x_i) \Bigr\}
\]
for $n\geq 2$, and define $\partial_n ^{\rm BR}=0$ for $n \leq 1$.

For example, we have
\begin{align*}
\partial_4^{\rm BR} (y, x_1, x_2, x_3, x_4) &=+\{(y, x_2, x_3, x_4)-(y*x_1, x_2 \oline* x_1, x_3 \oline* x_1, x_4 \oline* x_1)\}\\
&\phantom{=}-\{(y, x_1, x_3, x_4)-(y*x_2, x_1 \uline* x_2, x_3 \oline* x_2, x_4 \oline* x_2)\}\\
&\phantom{=} +\{(y, x_1, x_2, x_4)-(y*x_3, x_1 \uline* x_3, x_2 \uline* x_3, x_4 \oline* x_3)\}\\
&\phantom{=} -\{(y, x_1, x_2, x_3)-(y*x_4, x_1 \uline* x_4, x_2 \uline* x_4, x_3 \uline* x_4)\}.
\end{align*}

\begin{proposition}[cf.~\cite{JSC_ME_MS2004}]
$C_*^{\rm BR}(X; \mathbb Z)_Y:=(C_n^{\rm BR}(X; \mathbb Z)_Y,\partial_n^{\rm BR})_{n \in \mathbb Z}$ is a chain complex.
\end{proposition}

The \textit{birack chain complex} of $X$ with an $X$-set $Y$ is the chain complex $C_*^{\rm BR}(X; \mathbb Z)_Y$.  It determines  
the \textit{birack homology group} $H_n^{\rm BR}(X; \mathbb Z)_Y$ and the cohomology group $H^n_{\rm BR}(X; \mathbb Z)_Y$. 
For an abelian group $A$, the chain and cochain complexes $C_*^{\rm BR}(X; A)_Y$ and $C^*_{\rm BR}(X; A)_Y$ are also defined in the ordinary way.
We denote by $H_n^{\rm BR}(M; A)_Y$ and $H^n_{\rm BR}(X; A)_Y$ its homology and cohomology groups. 
A cocycle of $C^*_{\rm BR}(X;A)_Y$ is called a {\it birack cocycle} with $Y$ of the biquandle $X$ with $A$.
We omit $Y$, that is, we write $C_*^{\rm BR}(X; \mathbb Z):=C_*^{\rm BR}(X; \mathbb Z)_Y$ if $Y$ is the trivial $X$-set.
\vskip 1ex

\subsection{The (co)homology of multiple conjugation biquandles}\label{subsec_homology_MCB}
We review chain complexes  $P_*(X; \mathbb Z)_Y$, $D_*(X; \mathbb Z)_Y$ and $C_*(X; \mathbb Z)_Y$ defined in \cite{IIKKMO2} 
for multiple conjugation biquandles $X$ with $X$-set $Y$.  Refer to \cite{IIKKMO2} for a geometric interpretation of the chain complexes and an application in knot theory. These chain complexes are an analogy of chain complexes defined  in \cite{CarterIshiiSaitoTanaka16} for multiple conjugation quandles.

Let $X=\bigsqcup_{\lambda\in\Lambda}G_\lambda$ be a multiple conjugation biquandle and $Y$ an $X$-set of the MCB $X$. For $n \geq 0$, let $P_n(X; \mathbb Z)_Y$ denote the free abelian group generated by the elements
\[
\langle y \rangle \langle \bm{x}_1 \rangle \cdots \langle \bm{x}_k \rangle \in \bigcup_{n_1+\cdots+n_k=n}  Y \times \prod_{i=1}^k \bigcup_{\lambda \in \Lambda} G_\lambda^{n_i},
\]
where $\langle \boldsymbol{x}_i \rangle $ stands for an element $\langle x_{i1}, \ldots, x_{in_i} \rangle \in \bigcup_{\lambda \in \Lambda}G_\lambda^{n_i}$.
For $n < 0$, let $P_n(X; \mathbb Z)_Y=0$.
      We write
\begin{align*}
\langle \bm x_i  \oline* x \rangle &:= \langle x_{i1} \oline* x, x_{i2} \oline* x, \ldots, x_{in_i} \oline* x \rangle \text{ and}  \\
\langle \bm x_i  \uline* x \rangle &:= \langle x_{i1} \uline* x, x_{i2} \uline* x, \ldots, x_{in_i} \uline* x \rangle
\end{align*}
for an element $x \in X$ and an element $\langle \bm{x}_i \rangle = \langle x_{i1}, \ldots, x_{in_i} \rangle \in \bigcup_{\lambda \in \Lambda} G_\lambda^{n_i}$.

Define a boundary map $\partial_n:P_n(X; \mathbb Z)_Y\to P_{n-1}(X; \mathbb Z)_Y$ by
\begin{align*}
&
\partial_{n} ( \langle y \rangle \langle \bm{x}_1 \rangle  \cdots \langle \bm{x}_k \rangle ):=
\sum_{i=1}^{k} (-1)^{n_1+\cdots+n_{i-1} } \Biggl\{ \\
&
\hspace{2em}\langle y * x_{i1} \rangle
\langle \bm{x}_1 \uline* x_{i1} \rangle
\cdots
\langle \bm{x}_{i-1} \uline* x_{i1} \rangle
\langle \widetilde{\bm{x}}_i \oline* x_{i1}\rangle
\langle \bm{x}_{i+1}  \oline* x_{i1} \rangle
\cdots
\langle \bm{x}_k  \oline* x_{i1} \rangle \\
&
\hspace{2em}
+\sum_{j=1}^{n_i}
(-1)^{j}
\langle y \rangle
\langle \bm{x}_1 \rangle
\cdots
\langle \bm{x}_{i-1}  \rangle
\langle x_{i1}, \ldots, x_{i(j-1)}, x_{i(j+1)}, \ldots,x_{in_i} \rangle
\langle \bm{x}_{i+1} \rangle
\cdots
\langle \bm{x}_{k} \rangle
\Biggr \}
\end{align*}
when $n=n_1+\cdots+n_k \ge 1$, and define $\partial_n=0$ if $n<0$,
where 
\[
\langle \widetilde{\bm{x}}_i \oline* x_{i1}\rangle
:=
\langle 
x_{i1}^{-1} x_{i2} \oline* x_{i1}, 
x_{i1}^{-1} x_{i3} \oline* x_{i1}, \ldots, 
x_{i1}^{-1} x_{i(n_i-1)} \oline* x_{i1}, 
x_{i1}^{-1} x_{in_i} \oline* x_{i1} \rangle. 
\]

For example, the boundary maps are computed as follows.
\begin{align*}
&\partial_3(\langle y\rangle\langle x_1\rangle\langle x_2\rangle\langle x_3\rangle)\\
&=\langle y*x_1\rangle\langle x_2\oline{*}x_1\rangle\langle x_3\oline{*}x_1\rangle
-\langle y\rangle\langle x_2\rangle\langle x_3\rangle \\
&\hspace{5mm}
-\langle y*x_2\rangle\langle x_1\uline{*}x_2\rangle\langle x_3\oline{*}x_2\rangle
+\langle y\rangle\langle x_1\rangle\langle x_3\rangle \\
&\hspace{5mm}
+\langle y*x_3\rangle\langle x_1\uline{*}x_3\rangle\langle x_2\uline{*}x_3\rangle
-\langle y\rangle\langle x_1\rangle\langle x_2\rangle, \\
&\partial_3(\langle y\rangle\langle x_1\rangle\langle x_2,x_3\rangle)\\
&=\langle y*x_1\rangle\langle x_2\oline{*}x_1,x_3\oline{*}x_1\rangle
-\langle y\rangle\langle x_2,x_3\rangle  -\langle y*x_2\rangle\langle x_1\uline{*}x_2\rangle\langle x_2^{-1}x_3\oline{*}x_2\rangle \\
&\hspace{5mm}
+\langle y\rangle\langle x_1\rangle\langle x_3\rangle
-\langle y\rangle\langle x_1\rangle\langle x_2\rangle, \\
&\partial_3(\langle y\rangle\langle x_1,x_2\rangle\langle x_3\rangle)\\
&=\langle y*x_1\rangle\langle x_1^{-1}x_2\oline{*}x_1\rangle\langle x_3\oline{*}x_1\rangle
-\langle y\rangle\langle x_2\rangle\langle x_3\rangle \\
&\hspace{5mm}+\langle y\rangle\langle x_1\rangle\langle x_3\rangle
+\langle y*x_3\rangle\langle x_1\uline{*}x_3,x_2\uline{*}x_3\rangle
-\langle y\rangle\langle x_1,x_2\rangle, \\
&\partial_3(\langle y\rangle\langle x_1,x_2,x_3\rangle)\\
&=\langle y*x_1\rangle\langle x_1^{-1}x_2\oline{*}x_1,x_1^{-1}x_3\oline{*}x_1\rangle 
-\langle y\rangle\langle x_2,x_3\rangle
+\langle y\rangle\langle x_1,x_3\rangle
-\langle y\rangle\langle x_1,x_2\rangle.
\end{align*}

\begin{remark}
The notations of this paper are different from that of the paper \cite{IIKKMO2}, where 
we used more notations in order to prove propositions clearly. 
\end{remark}

\begin{proposition}[\cite{IIKKMO2}]
$P_*(X; \mathbb Z)_Y:=(P_n(X ; \mathbb Z)_Y,\partial_n)_{n \in \mathbb Z}$ is a chain complex.
\end{proposition}
For positive integers $s$ and $t$, we put 
$$M(s,t):=\{ \mu:\{1, \ldots, s \} \to \{ 1, \ldots, t \} ~|~ i < j \Rightarrow \mu(i) < \mu(j) \}.$$  
For any map $\mu \in M(s,t)$ and any integer $j \in \mathbb{Z}$, we define $\lfloor j;\mu\rfloor\in\{0, 1,\ldots, s\}$ as 
$\lfloor j;\mu\rfloor := \max \{\ell~|~\mu(\ell) \le j\}$.

For any $\lambda \in \Lambda$ and for elements $a_1, \ldots, a_s, b_1, \ldots, b_t \in G_\lambda$, put $\langle\boldsymbol{a}\rangle:=\langle a_1,\ldots,a_s\rangle$ and $\langle\boldsymbol{b}\rangle:=\langle b_1,\ldots,b_t\rangle$. We set
\begin{align*}
\big \langle \langle \boldsymbol{a} \rangle \langle \boldsymbol{b} \rangle \big \rangle_{\mu} &:=(-1)^{\sum_{k=1}^s(\mu(k)-k)}
\big \langle  a_{\lfloor 1;\mu\rfloor}b_{1-\lfloor 1; \mu\rfloor}, \ldots, a_{\lfloor s+t; \mu \rfloor}b_{s+t-\lfloor s+t;\mu \rfloor} \big \rangle \text{ and}\\
\big\langle \langle \boldsymbol{a}\rangle\langle\boldsymbol{b} \rangle \big\rangle &:= \displaystyle \sum_{\mu \in M(s, s+t)} \big \langle \langle \boldsymbol{a} \rangle \langle \boldsymbol{b} \rangle \big \rangle_{\mu},
\end{align*}
where $a_0=b_0=e_{\lambda}$.
For example, we have
\begin{align*}
\langle\langle a_1,a_2\rangle\langle b_1,b_2\rangle\rangle
&=\langle a_1,a_2,a_2b_1,a_2b_2\rangle
-\langle a_1,a_1b_1,a_2b_1,a_2b_2\rangle \\
&\hspace{1em}
+\langle b_1,a_1b_1,a_2b_1,a_2b_2\rangle
+\langle a_1,a_1b_1,a_1b_2,a_2b_2\rangle \\
&\hspace{1em}
-\langle b_1,a_1b_1,a_1b_2,a_2b_2\rangle
+\langle b_1,b_2,a_1b_2,a_2b_2\rangle.
\end{align*}
for elements $a_1, a_2, b_1, b_2 \in G_{\lambda}$.

\begin{definition}
Let $D_n(X;\mathbb Z)_Y$ be the submodule of $P_n(X ; \mathbb Z )_Y$ generated by elements of the form $\langle y\rangle\langle\boldsymbol{a}_1\rangle\cdots
\langle\boldsymbol{a}\rangle\langle\boldsymbol{b}\rangle
\cdots\langle\boldsymbol{a}_k\rangle
-\langle y\rangle\langle\boldsymbol{a}_1\rangle\cdots
\langle\langle\boldsymbol{a}\rangle\langle\boldsymbol{b}\rangle\rangle
\cdots\langle\boldsymbol{a}_k\rangle$.
When $n\leq1$, we set $D_n(X ; \mathbb Z )_Y=0$.
\end{definition}

For example, the submodule $D_2(X ; \mathbb Z )_Y$ is generated by the elements of the form
\[ \langle y\rangle\langle a\rangle\langle b\rangle
-\langle y\rangle\langle a,ab\rangle
+\langle y\rangle\langle b,ab\rangle, \]
for $y\in Y$, $a,b\in G_\lambda$.

\begin{proposition}[\cite{IIKKMO2}]
$D_*(X; \mathbb Z)_Y:=(D_n(X; \mathbb Z)_Y,\partial_n)_{n \in \mathbb Z}$ is a subcomplex of $P_*(X; \mathbb Z)_Y$.
\end{proposition}

We have a chain complex $C_*(X;\mathbb Z)_Y:=P_*(X;\mathbb Z)_Y/D_*(X;\mathbb Z)_Y$.  
The homology group 
$H_n(X; \mathbb Z)_Y$ and the cohomology group $H^n(X; \mathbb Z)_Y$ 
of the multiple conjugation biquandle $X$ with $X$-set $Y$ are the homology and cohomology groups of the chain complex $C_*(X;\mathbb Z)_Y$. 

For an abelian group $A$, the chain and cochain complexes $C_*(X; A)_Y$ and $C^*(X; A)_Y$ are also defined in the ordinary way.
We denote by $H_n(X; A)_Y$ and $H^n(X; A)_Y$ its $n$-th homology and cohomology groups respectively.
A cocycle of $C^*(X;A)_Y$ is called a {\it cocycle} with $Y$ of the multiple conjugation biquandle $X$ with  $A$. We omit $Y$, that is, we write $C_*(X; \mathbb Z):=C_*(X; \mathbb Z)_Y$ if $Y$ is the trivial $X$-set.
\vskip 1ex

\subsection{The (co)homology of groups}\label{subsec_homology_group}
In this subsection, we review chain complexes for groups.

We review the group (co)homology. Let $G$ be a group and let $M$ be a right $\mathbb Z[G]$-module, where $\mathbb Z[G]$ is the group ring of $G$ over the ring $\mathbb Z$. Let $C^{\rm gp}_n( M ; \mathbb Z)$  be a  free abelian group generated by any element $(m_1, \ldots, m_n) \in M^n$ if $n\geq 1$, and $C^{\rm gp}_n(M ; \mathbb Z)=0$ otherwise.

We define the boundary map $\partial_n ^{\rm gp}: C^{\rm gp}_n(M;\mathbb Z) \to C^{\rm gp}_{n-1}(M;\mathbb Z)$ by 
\begin{align*}
\partial_n ^{\rm gp}(m_1, \ldots , m_n) &:= (m_2, \ldots, m_n)\\
&\phantom{=}+\displaystyle \sum_{i=1}^{n-1}  (-1)^i (m_1, \ldots , m_{i-1},m_i+m_{i+1},m_{i+2}, \ldots, m_n)\\
&\phantom{=}+(-1)^n(m_1, \ldots, m_{n-1})
\end{align*}
if $n\geq 2$, and $\partial_n ^{\rm gp}=0$ otherwise.
For example,
\begin{align*}
\partial_4 ^{\rm gp}(m_1, m_2, m_3, m_4) &= (m_2, m_3, m_4)-(m_1+m_2, m_3, m_4)+(m_1, m_2+m_3, m_4)\\
&\phantom{=}-(m_1, m_2, m_3+m_4)+(m_1, m_2, m_3).
\end{align*}

Let 
$$C^{\rm gp}_n(M;\mathbb Z)_G:=C^{\rm gp}_n(M;\mathbb Z) \otimes_{\mathbb Z[G]} {\mathbb Z},$$ 
where $G$ acts on $\mathbb Z$ trivially from the left. 
In other words, 
under the diagonal action of $G$ on  $C^{\rm gp}_n(M;\mathbb Z)$, 
$C^{\rm gp}_n(M;\mathbb Z)_G$ is the $G$-coinvariant part. We note that $(m_1g, \ldots, m_ng)=(m_1, \ldots, m_n)$ in $C^{\rm gp}_n(M;\mathbb Z)_G$ for any $(m_1, \ldots, m_n) \in M^n$ and any $g \in G$. We have the induced boundary map $\partial_n ^{\rm gp}: C^{\rm gp}_n(M;\mathbb Z)_G \to C^{\rm gp}_{n-1}(M;\mathbb Z)_G$.

\begin{proposition}[\cite{EM1947}]
$C_*^{\rm gp}(M;\mathbb Z)_G=(C_n^{\rm gp}(M;\mathbb Z)_G, \partial_n^{\rm gp})_{n \in \mathbb Z}$ is a chain complex. 
\end{proposition}

The chain complex determines the homology group $H_n^{\rm gp}(M; \mathbb Z)_G$ and the cohomology group $H^n_{\rm gp}(M; \mathbb Z)^G$.
 For an abelian group $A$, the chain and cochain complexes $C_*^{\rm gp}(M; A)_G$ and $C^*_{\rm gp}(M; A)^G$ are also defined in the ordinary way.
We denote by $H_n^{\rm gp}(M; A)_G$ and $H^n_{\rm gp}(M; A)^G$ its homology and cohomology groups respectively.

Let $D^{\rm gp}_n(M;\mathbb Z)_G$ be  the submodule of $C^{\rm gp}_n(M;\mathbb Z)_G$ generated by the elements of the following set
\[
\displaystyle \bigcup_{ i=1}^{ n-1} \Big\{ (m_1, \ldots, m_{i-1},0, m_{i+1},\ldots m_n ) ~\Big|~m_1, \ldots, m_n \in M \Big\}
\]
for $n \ge 2$.
Define $D^{\rm gp}_n(M;\mathbb Z)_G=0$ for  $n \le 1$.

\begin{proposition}[\cite{EM1947}]
$D^{\rm gp}_*(M;\mathbb Z)_G:=(D^{\rm gp}_n(M;\mathbb Z)_G, \partial_n^{\rm gp})_{n \in \mathbb Z}$ is a subcomplex of $C^{\rm gp}_*(M;\mathbb Z)_G$.
\end{proposition}

The chain complex $C^{\rm norgp}_*(M;\mathbb Z)_G:=C^{\rm gp}_*(M;\mathbb Z)_G/D^{\rm gp}_*(M;\mathbb Z)_G$ determines the homology group $H^{\rm norgp}_n (M; \mathbb Z)_G$ called the {\it normalized group homology}.
In the ordinary way, for an abelian group $A$,  the homology group $H^{\rm norgp}_n (M; A)_G$ and the cohomology group $H_{\rm norgp}^n (M; A)_G$ are defined. 	
\vspace{1ex}

\section{Cocycles of $G$-Alexander biquandles with trivial $X$-set}\label{app:11}

Throughout this section, let $X= M \times G$ be the $G$-Alexander  biquandle of $(M, \varphi)$, see Definition \ref{def_Gfamily_of_Alexander_quandle}. 
We suppose that $X$ has trivial $X$-set.  

We define a chain complex $C^{\rm BR_U}_\ast (X;\mathbb Z)$ and chain maps $\gamma$ and $\psi$. 
It turns out that the chain complex  $C^{\rm BR_U}_\ast (X;\mathbb Z)$ is isomorphic to the birack chain complex  $C^{\rm BR}_\ast (X;\mathbb Z)$, and hence the homology group $H^{\rm BR_U}_n (X;\mathbb Z)$ is isomorphic to  $H^{\rm BR}_n (X;\mathbb Z)$.  Through this chain complex, we associate cocycles of the chain complex of group $G$ to cocycles of the birack chain complex. 

Our goal in this section is to give Theorems~\ref{birack_cocycle_main} and~\ref{birack cocycle explicit}.

\subsection{The chain complex $C^{\rm BR_U}_\ast (X;\mathbb Z)$ and chain maps $\gamma$ and $\psi$}\label{subsection:A1}

\subsubsection{The chain complex $C^{\rm BR_U}_\ast (X;\mathbb Z)$ and the homology group $H^{\rm BR_U}_{n}(X; \mathbb Z)$}

For $\bm{g}=(g_1,  \ldots , g_n) \in G^n$ and $\bm{m}=(m_1, \ldots, m_n) \in M^n$, we use the following notations:
{\normalsize
\begin{align*}
\bm{g}_{\{i\}}:=&(g_1,  \ldots, g_{i},g_{i+2},\ldots , g_n) \in G^{n-1},\\
\bm{g}_{\{ \triangleleft i\} }:=&(g_{i+1}^{-1} g_1 g_{i+1},  \ldots, g_{i+1}^{-1} g_{i} g_{i+1},g_{i+2},\ldots , g_n) \in G^{n-1},\\
\bm{m}_{\{i\}}:=&(m_1, \ldots, m_{i-1},m_{i} + m_{i+1},m_{i+2}, \ldots , m_n) \in M^{n-1} \text{ and}\\
\bm{m}_{\{ \triangleleft i\} }:=&(m_1 g_{i+1}, \ldots, m_{i-1} g_{i+1},m_{i} g_{i+1} + m_{i+1} \varphi(g_{i+1}),\\
&\hspace{12ex}m_{i+2} \varphi(g_{i+1}), \ldots , m_n \varphi(g_{i+1})) \in M^{n-1},
\end{align*}
}
where $n$ and $i$ are integers such that $n \ge 2$ and $0 \le i \le n-1$.

Let $C^{\rm BR_U}_n(X; \mathbb Z)$ be the free abelian group generated by the elements $(\bm{g};\bm{m})=(g_1, \ldots, g_n; m_1, \ldots, m_n)\in G^{n} \times M^{n}$ if $n\geq 1$, and $C^{\rm BR_U}_n(X; \mathbb Z)=0$ otherwise. 
Define a boundary map $\partial_n^{\rm BR_U}: C^{\rm BR_U}_n(X; \mathbb Z) \to C^{\rm BR_U}_{n-1}(X; \mathbb Z)$ by 
\[
\begin{array}{ll}
\partial_n^{\rm BR_U} (\bm{g};\bm{m})
=\displaystyle \sum_{i=0}^{n-1} (-1)^i  \left\{ (\bm{g}_{\{i\}};\bm{m}_{\{i\}})- ( \bm{g}_{\{ \triangleleft i \}};\bm{m}_{\{ \triangleleft i \}}) \right\}
\end{array}
\]
for $n\geq 2$,  and define $\partial_n^{\rm BR_U}=0$ for $n \leq 1$. For example, we have
\begin{align*}
&\partial_3^{\rm BR_U} (g_1, g_2, g_3; m_1, m_2, m_3)\\
&=\bigl\{ (g_2, g_3; m_2, m_3)-(g_2, g_3; m_2\varphi(g_1), m_3\varphi(g_1))\bigr\}\\
&\phantom{=}-\bigl\{(g_1, g_3; m_1+m_2, m_3)-(g_2^{-1}g_1g_2, g_3; m_1g_2+m_2\varphi(g_2), m_3\varphi(g_2)) \bigr\}\\
&\phantom{=}+\bigl\{ (g_1, g_2; m_1, m_2+m_3) - (g_3^{-1}g_1g_3, g_3^{-1}g_2g_3; m_1g_3, m_2g_3+m_3\varphi(g_3)) \bigr\}.
\end{align*}
\begin{lemma}
$C^{\rm BR_U}_*(X; \mathbb Z):=(C^{\rm BR_U}_n(X; \mathbb Z), \partial_n^{\rm BR_U})_{n \in \mathbb Z}$ is a chain complex. 
\end{lemma}
\begin{proof}
We fix an integer $n$ with $n \ge 2$. Define  $\partial_{n}^{i}: C^{\rm BR_U}_n(X; \mathbb Z) \to C^{\rm BR_U}_{n-1}(X; \mathbb Z)$ by $\partial_{n}^{i} (\bm{g};\bm{m}):=(\bm{g}_{\{i\}};\bm{m}_{\{i\}})$  and ${\delta}_{n}^{i}: C^{\rm BR_U}_n(X; \mathbb Z) \to C^{\rm BR_U}_{n-1}(X; \mathbb Z)$ by 
$\delta_{n}^{i} (\bm{g}; \bm{m}):=(\bm{g}_{\{ \triangleleft i \}} ; \bm{m}_{\{ \triangleleft i \}})$ 
for any integer $i$ with $0\leq i \leq n-1$. Then, we have
\[
\partial_n^{\rm BR_U}=\displaystyle \sum_{i=0}^{n-1} (-1)^i (\partial_{n}^i - \delta_{n}^{i}).
\]
We assert that
\begin{align}
\delta_{n-1}^j \circ \delta_{n}^i &= \delta_{n-1}^{i} \circ \delta_{n}^{j+1}, \label{chain1}\\
\partial_{n-1}^j \circ \partial_{n}^i &= \partial_{n-1}^{i} \circ \partial_{n}^{j+1},\\
\partial_{n-1}^j \circ \delta_{n}^i &= \delta_{n-1}^{i} \circ \partial_{n}^{j+1} \text{ and}\\
\delta_{n-1}^j \circ \partial_{n}^i &= \partial_{n-1}^{i} \circ \delta_{n}^{j+1}.\label{chain2}
\end{align}
We show the equality \eqref{chain1}.
Since the other equalities are shown in a similar way,
we omit the proofs.

Let $i$ and $j$ be integers with $0\leq i \leq j\leq n-2$. We see $\bm{g}_{\{ \triangleleft i\}\{ \triangleleft j\}}=\bm{g}_{\{ \triangleleft (j+1)\}\{ \triangleleft i\}}$ as follows: 
\begin{align*}
\bm{g}_{\{ \triangleleft i\}\{ \triangleleft j\}}&=(g_{j+2}^{-1}g_{i+1}^{-1}g_1g_{i+1}g_{j+2},\ldots,g_{j+2}^{-1}g_{i+1}^{-1}g_ig_{i+1}g_{j+2},g_{j+2}^{-1}g_{i+2}g_{j+2},\ldots,\\
&\phantom{=(}g_{j+2}^{-1}g_{j+1}g_{j+2},g_{j+3},\ldots,g_n)\\
&=((g_{j+2}^{-1}g_{i+1}g_{j+2})^{-1}g_{j+2}^{-1}g_1g_{j+2}(g_{j+2}^{-1}g_{i+1}g_{j+2}),\ldots,\\
&\phantom{=(}(g_{j+2}^{-1}g_{i+1}g_{j+2})^{-1}g_{j+2}^{-1}g_{i}g_{j+2}(g_{j+2}^{-1}g_{i+1}g_{j+2}),g_{j+2}^{-1}g_{i+2}g_{j+2},\ldots,\\
&\phantom{=(}g_{j+2}^{-1}g_{j+1}g_{j+2},g_{j+3},\ldots,g_n)\\
&=\bm{g}_{\{ \triangleleft (j+1)\}\{ \triangleleft i\}}.
\end{align*}
We see
$\delta_{n-1}^j \circ \delta_{n}^i = \delta_{n-1}^{i} \circ \delta_{n}^{j+1}$ for $0<i<j\leq n-2$ as follows: 
\begin{align*}
&\delta_{n-1}^{j} \circ \delta_{n}^{i}  (\bm{g}; \bm{m}) \\
&= (\bm{g}_{\{ \triangleleft i\}\{ \triangleleft j\}}; m_1g_{i+1}g_{j+2},\ldots,m_{i-1}g_{i+1}g_{j+2}, m_ig_{i+1}g_{j+2}+m_{i+1}\varphi(g_{i+1})g_{j+2},\\
&\phantom{=(} m_{i+2}\varphi(g_{i+1})g_{j+2},\ldots,m_{j}\varphi(g_{i+1})g_{j+2}, m_{j+1}\varphi(g_{i+1})g_{j+2}+m_{j+2}\varphi(g_{i+1})\varphi(g_{j+2}),\\
&\phantom{=(} m_{j+3}\varphi(g_{i+1})\varphi(g_{j+2}),\ldots, m_{n}\varphi(g_{i+1})\varphi(g_{j+2}))\\
&= (\bm{g}_{\{ \triangleleft (j+1)\}\{ \triangleleft i\}}; m_1g_{j+2}(g_{j+2}^{-1}g_{i+1}g_{j+2}),\ldots,m_{i-1}g_{j+2}(g_{j+2}^{-1}g_{i+1}g_{j+2}),\\
&\phantom{=(} m_ig_{j+2}(g_{j+2}^{-1}g_{i+1}g_{j+2})+m_{i+1}g_{j+2}\varphi(g_{j+2}^{-1}g_{i+1}g_{j+2}),m_{i+2}g_{j+2}\varphi(g_{j+2}^{-1}g_{i+1}g_{j+2}),\ldots,\\
&\phantom{=(} m_{j}g_{j+2}\varphi(g_{j+2}^{-1}g_{i+1}g_{j+2}),m_{j+1}g_{j+2}\varphi(g_{j+2}^{-1}g_{i+1}g_{j+2})+m_{j+2}\varphi(g_{j+2})\varphi(g_{j+2}^{-1}g_{i+1}g_{j+2}),\\
&\phantom{=(} m_{j+3}\varphi(g_{j+2})\varphi(g_{j+2}^{-1}g_{i+1}g_{j+2}),\ldots, m_{n}\varphi(g_{j+2})\varphi(g_{j+2}^{-1}g_{i+1}g_{j+2}))\\
&=\delta_{n-1}^{i} \circ \delta_{n}^{j+1}  (\bm{g}; \bm{m}).
\end{align*}
For the cases that  $i=j$ or $i=0$, we can show that the equality \eqref{chain1} holds.

Hence, we have the equality \eqref{chain1} for any $i,j$ with $0\leq i \leq j\leq n-2$.

By using the equalities \eqref{chain1}--\eqref{chain2}, we see that
$\partial_{n-1}^{\rm BR_U}\circ\partial_n^{\rm BR_U}=0$.
\end{proof}

The $n$-th homology and cohomology groups of the chain complex $C_*^{\rm BR_U} (X; \mathbb Z)$ are denoted by  
$H_n^{\rm BR_U}(X; \mathbb Z)$ and $H^n_{\rm BR_U}(X; \mathbb Z)$ respectively.



\subsubsection{The chain map $\gamma: C_*^{\rm BR} (X; \mathbb Z) \to C_{*}^{\rm BR_U} (X; \mathbb Z)$}

For $n\geq 1$, define a homomorphism $\gamma_n: C_n^{\rm BR} (X; \mathbb Z) \to C_{n}^{\rm BR_U} (X; \mathbb Z)$ by 
\[
\gamma_n((m_1, g_1), \ldots , (m_n, g_n)) := (\bm{g};m_1', \ldots, m_{n-1}', m_n),
\]
where we write $m'_i:=m_i - m_{i+1}$ and $\bm{g} = (g_1, \dots, g_n)$. 
For $n < 1$, define $\gamma_n = 0$.

\begin{lemma}\label{3lem-gamchain} 
\begin{itemize}
\item[(1)] 
The map $\gamma$ is a chain map, that is, it holds that
\[ \gamma_{n-1} \circ \partial_n^{\rm BR}= \partial_n^{\rm BR_U} \circ \gamma_{n}. \]

\item[(2)] For any integer $n$, the map $\gamma_n$ is an isomorphism.  
\end{itemize}
\end{lemma}

\begin{proof} 

The assertion that $\gamma$ is a chain map is easily verified by a direct calculation, which is left to the reader.  
Define a homomorphism $r_n:  C_{n}^{\rm BR_U} (X; \mathbb Z) \to C_n^{\rm BR} (X; \mathbb Z) $ by
\[
r_n(\bm{g};\bm{m}):=
\biggl(
\Bigl(\sum_{i=1}^{n}m_i,~g_1 \Bigr),~\Bigl( \sum_{i=2}^{n}m_i,~g_2 \Bigr),\ldots,~\Bigl(\sum_{i=n}^{n}m_i,~g_n \Bigr)
\biggr).
\]
We can easily see that $r_n$ is the inverse map of $\gamma_n$, and thus, $\gamma_n$ is a bijection.
\end{proof}

\subsubsection{The chain map $\psi : C_{*}^{\rm BR_U} (X; \mathbb Z) \to C_{*}^{\rm gp} (M; \mathbb Z)_G$}

For $n\geq 1$, define $\psi_n : C_{n}^{\rm BR_U} (X; \mathbb Z) \to C_{n}^{\rm gp} (M; \mathbb Z)_G$ by 
\[
\psi_n (\bm{g};\bm{m}):=
\displaystyle \sum_{\bm{k} =(k_1, \ldots, k_n) \in \mathcal K_n} (-1)^{| \bm{k} |} (m_1  \bm{g}_{\bm{k}_1}, \ m_2  \bm{g}_{\bm{k}_2}, \ \ldots, \ m_n  \bm{g}_{\bm{k}_n}),
\]
where $\mathcal{K}_n :=\{\bm{k} = (k_1, \ldots , k_n) \in \{0,1\}^n ~|~ k_1=0\}$
and 
\begin{align*}
|\bm{k}| &:= k_1 + k_2 + \cdots  +k_n \text{ and} \\
\bm{g}_{\bm{k}_i} &:= \varphi(g_1^{k_1} g_2^{k_2} \cdots g_{i}^{k_i}) g_{i+1}^{k_{i+1}}g_{i+2}^{k_{i+2}} \cdots g_{n}^{k_{n}} \in G
\end{align*}
for an element $\bm{k}=(k_1, \ldots, k_n) \in \mathcal K_n$ and an integer $i$ with $1 \le i \le n$.
For $n <1$, define $\psi_n=0$.

For example, we have 
\begin{align*}
\psi_3(\bm{g};\bm{m}) &= \phantom{+} (-1)^{0+1+1}(m_1\varphi(g_1^0)g_2^1g_3^1, m_2\varphi(g_1^0g_2^1)g_3^1, m_3\varphi(g_1^0g_2^1g_3^1) )\\
&\phantom{=}+ (-1)^{0+1+0}(m_1\varphi(g_1^0)g_2^1g_3^0, m_2\varphi(g_1^0g_2^1)g_3^0, m_3\varphi(g_1^0g_2^1g_3^0) )\\
&\phantom{=}+ (-1)^{0+0+1}(m_1\varphi(g_1^0)g_2^0g_3^1, m_2\varphi(g_1^0g_2^0)g_3^1, m_3\varphi(g_1^0g_2^0g_3^1) )\\
&\phantom{=}+ (-1)^{0+0+0}(m_1\varphi(g_1^0)g_2^0g_3^0, m_2\varphi(g_1^0g_2^0)g_3^0, m_3\varphi(g_1^0g_2^0g_3^0) )\\
&= \phantom{+} (m_1g_2g_3, m_2\varphi(g_2)g_3, m_3\varphi(g_2g_3) )\\
&\phantom{=}-(m_1g_2, m_2\varphi(g_2), m_3\varphi(g_2) )\\
&\phantom{=}-(m_1g_3, m_2g_3, m_3\varphi(g_3) )\\
&\phantom{=}+(m_1, m_2, m_3 ).
\end{align*}
\begin{lemma}
The map $\psi$ is a chain map, that is, it holds that
\[ \psi_{n-1} \circ \partial_n^{\rm BR_U}= \partial_n^{\rm gp} \circ \psi_{n}.\]
\end{lemma}

\begin{proof}
We fix an integer $n \ge 2$, and  
fix $\bm{g} =(g_1, \dots, g_n)$ and $\bm{m} =(m_1, \dots,m_n)$. 

For any integer $i$ with $2\leq i \leq n$, let
\begin{align*}
\mathcal{K}_{i}^{0} &:= \{\bm{k} = (k_1, \ldots , k_n) \in \{0,1\}^n ~|~ k_1=k_i=0\} \text{ and} \\
\mathcal{K}_{i}^{1} &:= \{\bm{k} = (k_1, \ldots , k_n) \in \{0,1\}^n ~|~ k_1=0, k_i=1\}.
\end{align*}
For any $i$ with $1 \le i \le n-1$, put
\begin{align*}
A_i &:=( m_1 \bm{g}_{\bm{k}_1},\ldots , m_{i-1} \bm{g}_{\bm{k}_{i-1}}, m_{i}  \bm{g}_{\bm{k}_{i}} + m_{i+1} \bm{g}_{\bm{k}_{i+1}}, m_{i+2} \bm{g}_{\bm{k}_{i+2}} , \ldots ,m_n   \bm{g}_{\bm{k}_{n}}\big),\\
L &:=\displaystyle \sum_{\bm{k} \in \mathcal{K}_n} (-1)^{|\bm{k}|} \displaystyle \sum_{i=1}^{n-1} (-1)^i A_i.
\end{align*}
We show that
$
\psi_{n-1}\circ \partial_n^{\rm BR_U} (\bm{g}; \bm{m})=L=\partial_{n}^{\rm gp}\circ\psi_{n}(\bm{g};\bm{m})$.

Firstly, we show $\psi_{n-1}\circ \partial_n^{\rm BR_U} (\bm{g}; \bm{m})=L$.
We have
\begin{align}
\psi_{n-1}\circ \partial_n^{\rm BR_U} (\bm{g}; \bm{m}) &=\psi_{n-1} \big((\bm{g}_{\{0\}};\bm{m}_{\{0\}})- ( \bm{g}_{\{ \triangleleft 0 \}} ; \bm{m}_{\{ \triangleleft 0 \}})\big) \nonumber \\
&\phantom{=} + \displaystyle \sum_{i=1}^{n-1} (-1)^i \psi_{n-1} \big((\bm{g}_{\{i\}};\bm{m}_{\{i\}})- ( \bm{g}_{\{ \triangleleft i \}} ; \bm{m}_{\{ \triangleleft i \}})\big). \label{cad} 
\end{align}
Since $(m_1g, \ldots, m_ng)=(m_1, \ldots, m_n)$ in $C^{\rm gp}_n(M;\mathbb Z)_G$ for any $(m_1, \ldots, m_n) \in M^n$ and any $g \in G$, we see that
\begin{align}
&\psi_{n-1}(\bm{g}_{\{0\}};\bm{m}_{\{0\}})- \psi_{n-1}(\bm{g}_{\{ \triangleleft 0 \}} ; \bm{m}_{\{ \triangleleft 0 \}})  \nonumber \\
&= \displaystyle \sum_{\bm{k} \in \mathcal{K}_{2}^0}(-1)^{|\bm{k}|} \big\{ (m_2 \bm{g}_{\bm{k}_2} , \ldots, m_n \bm{g}_{\bm{k}_n}) - (m_2 \varphi(g_1) \bm{g}_{\bm{k}_2} , \ldots, m_n \varphi(g_1)  \bm{g}_{\bm{k}_n} )\big\}  \nonumber \\
&= \displaystyle \sum_{\bm{k} \in \mathcal{K}_{2}^0}(-1)^{|\bm{k}|} \big\{ (m_2 \bm{g}_{\bm{k}_2} , \ldots, m_n \bm{g}_{\bm{k}_n}) - (m_2 \bm{g}_{\bm{k}_2} \varphi(g_1), \ldots, m_n \bm{g}_{\bm{k}_n} \varphi(g_1))\big\} \nonumber \\
&=0 \label{equality02}.
\end{align}
For any integer $i$ with $1 \le i \le n-1$, we have
\begin{align}
&\psi_{n-1}(\bm{g}_{\{i\}} ; \bm{m}_{\{i\}}) \nonumber \\
&= \displaystyle \sum_{\bm{k} \in \mathcal{K}_{i+1}^0} (-1)^{|\bm{k}|}  \Bigl( \nonumber \\
&\begin{array}{l@{}r@{}l@{}}
\phantom{=} & m_1 &\varphi(g_1^{k_1}) g_2^{k_2} \cdots  g_{i}^{k_{i}} g_{i+1}^0 g_{i+2}^{k_{i+2}} \cdots g_n^{k_n},\ldots, \\
\phantom{=} & m_{i-1} &\varphi(g_1^{k_1} \cdots g_{i-1}^{k_{i-1}}) g_{i}^{k_{i}} g_{i+1}^0 g_{i+2}^{k_{i+2}} \cdots g_n^{k_n}, \\
\phantom{=} & m_{i}  &\varphi(g_1^{k_1} \cdots g_{i}^{k_{i}} )  g_{i+1}^0 g_{i+2}^{k_{i+2}} \cdots g_n^{k_n}+m_{i+1} \varphi(g_1^{k_1}  \cdots  g_{i}^{k_{i}}   g_{i+1}^0 ) g_{i+2}^{k_{i+2}} \cdots g_n^{k_n},\\
\phantom{=} & m_{i+2} &\varphi(g_1^{k_1}  \cdots g_{i}^{k_{i}}   g_{i+1}^0 g_{i+2}^{k_{i+2}} ) g_{i+3}^{k_{i+3}} \cdots g_n^{k_n},\ldots, \\
\phantom{=} & m_{n}  &\varphi(g_1^{k_1} g_2^{k_2} \cdots g_{i}^{k_{i}}   g_{i+1}^0 g_{i+2}^{k_{i+2}}  \cdots g_n^{k_n}) \Bigl)
\end{array} \nonumber \\
&= \displaystyle \sum_{\bm{k} \in \mathcal{K}_{i+1}^0} (-1)^{|\bm{k}|} A_i. \label{equation02}
\end{align}
For any integer $i$ with $1 \le i \le n-1$, we also have
{\normalsize
\begin{align}
&\psi_{n-1}(\bm{g}_{\{ \triangleleft i \}} ; \bm{m}_{\{ \triangleleft i \}}) \nonumber \\
&= \displaystyle \sum_{\bm{k} \in \mathcal{K}_{i+1}^0} (-1)^{|\bm{k}|} \Bigl( \nonumber \\
&\begin{array}{l@{}r@{}l@{}}
\phantom{=} & m_1 g_{i+1} &\varphi \bigl( ( g_{i+1}^{-1} g_1 g_{i+1} )^{k_1} \bigr) ( g_{i+1}^{-1} g_2 g_{i+1} )^{k_2} \cdots ( g_{i+1}^{-1} g_{i} g_{i+1} )^{k_{i}} g_{i+1}^0 g_{i+2}^{k_{i+2}} \cdots g_n^{k_n},\ldots, \nonumber \\
\phantom{=} & m_{i-1} g_{i+1} &\varphi \bigl( ( g_{i+1}^{-1} g_1 g_{i+1} )^{k_1}  \cdots ( g_{i+1}^{-1} g_{i-1} g_{i+1} )^{k_{i-1}} \bigr) ( g_{i+1}^{-1} g_{i} g_{i+1} )^{k_{i}} g_{i+1}^0 g_{i+2}^{k_{i+2}} \cdots g_n^{k_n}, \nonumber \\
\phantom{=} & m_i g_{i+1} &\varphi \bigl( ( g_{i+1}^{-1} g_1 g_{i+1} )^{k_1}   \cdots ( g_{i+1}^{-1} g_{i} g_{i+1} )^{k_{i}} \bigr) g_{i+1}^0 g_{i+2}^{k_{i+2}} \cdots g_n^{k_n} \nonumber \\
\phantom{=} & +m_{i+1} \varphi(g_{i+1}) &\varphi \bigl( ( g_{i+1}^{-1} g_1 g_{i+1} )^{k_1}   \cdots ( g_{i+1}^{-1} g_{i} g_{i+1} )^{k_{i}} g_{i+1}^0 \bigr) g_{i+2}^{k_{i+2}} \cdots g_n^{k_n}, \nonumber \\
\phantom{=} & m_{i+2} \varphi(g_{i+1}) &\varphi \bigl( ( g_{i+1}^{-1} g_1 g_{i+1} )^{k_1}   \cdots  ( g_{i+1}^{-1} g_{i} g_{i+1} )^{k_{i}} g_{i+1}^0 g_{i+2}^{k_{i+2}} \bigr) g_{i+3}^{k_{i+3}} \cdots g_n^{k_n},\ldots, \\
\phantom{=} & m_n \varphi(g_{i+1}) &\varphi \bigl( ( g_{i+1}^{-1} g_1 g_{i+1} )^{k_1}  \cdots  ( g_{i+1}^{-1} g_{i} g_{i+1} )^{k_{i}} g_{i+1}^0 g_{i+2}^{k_{i+2}} \cdots g_n^{k_n} \bigr) \Bigl) \nonumber \\
\end{array}\\
&= \displaystyle \sum_{\bm{k} \in \mathcal{K}_{i+1}^0} (-1)^{|\bm{k}|} \Bigl( \nonumber \\
&\begin{array}{l@{}r@{}l@{}}
\phantom{=} & m_1 g_{i+1} &\varphi \bigl( g_1^{k_1} \bigr) g_{i+1}^{-1} g_2^{k_2} \cdots  g_{i}^{k_{i}} g_{i+1}^1 ~ g_{i+1}^0 g_{i+2}^{k_{i+2}} \cdots g_n^{k_n},\ldots, \\
\phantom{=} & m_{i-1}g_{i+1}  &\varphi \bigl( g_1^{k_1}  g_2^{k_2} \cdots  g_{i-1}^{k_{i-1}} \bigr)  g_{i+1}^{-1} g_{i}^{k_{i}} g_{i+1}^1 ~ g_{i+1}^0  g_{i+2}^{k_{i+2}} \cdots g_n^{k_n}, \\
\phantom{=} & m_i g_{i+1} &\varphi \bigl(  g_1^{k_1} g_2^{k_2} \cdots  g_{i}^{k_{i}} \bigr) g_{i+1}^0 g_{i+2}^{k_{i+2}} \cdots g_n^{k_n} \\
\phantom{=} & + m_{i+1} \varphi(g_{i+1}) &\varphi \bigl( g_1^{k_1} g_2^{k_2} \cdots  g_{i}^{k_{i}} g_{i+1}^0 \bigr) g_{i+2}^{k_{i+2}} \cdots g_n^{k_n}, \\
\phantom{=} & m_{i+2} \varphi(g_{i+1}) &\varphi \bigl( g_1^{k_1} g_2^{k_2} \cdots g_{i}^{k_{i}} g_{i+1}^0 g_{i+2}^{k_{i+2}} \bigr)g_{i+3}^{k_{i+3}} \cdots g_n^{k_n},\ldots, \\
\phantom{=} & m_n \varphi(g_{i+1})  &\varphi \bigl( g_1^{k_1} g_2^{k_2} \cdots   g_{i}^{k_{i}} g_{i+1}^0 g_{i+2}^{k_{i+2}} \cdots g_n^{k_n} \bigr) \Bigl)
\end{array} \nonumber \\
&= \displaystyle \sum_{\bm{k} \in \mathcal{K}_{i+1}^0} (-1)^{|\bm{k}|} \Bigl( \nonumber \\
&\begin{array}{l@{}r@{}l@{}}
\phantom{=} & m_1 &\varphi \bigl( g_1^{k_1} \bigr)  g_2^{k_2} \cdots  g_{i}^{k_{i}} g_{i+1}^1 g_{i+2}^{k_{i+2}} \cdots g_n^{k_n},\ldots, \\
\phantom{=} & m_{i-1} &\varphi \bigl( g_1^{k_1}  g_2^{k_2} \cdots  g_{i-1}^{k_{i-1}} \bigr)   g_{i}^{k_{i}} g_{i+1}^1  g_{i+2}^{k_{i+2}} \cdots g_n^{k_n}, \\
\phantom{=} & m_i   &\varphi \bigl(  g_1^{k_1} g_2^{k_2} \cdots   g_{i}^{k_{i}} \bigr) g_{i+1}^1 g_{i+2}^{k_{i+2}} \cdots g_n^{k_n} + m_{i+1}  \varphi \bigl( g_1^{k_1} g_2^{k_2} \cdots   g_{i}^{k_{i}} g_{i+1}^1 \bigr) g_{i+2}^{k_{i+2}} \cdots g_n^{k_n}, \\
\phantom{=} & m_{i+2}  &\varphi \bigl( g_1^{k_1} g_2^{k_2} \cdots  g_{i}^{k_{i}} g_{i+1}^1 g_{i+2}^{k_{i+2}} \bigr) g_{i+3}^{k_{i+3}} \cdots g_n^{k_n},\ldots, \\
\phantom{=} & m_n &\varphi \bigl( g_1^{k_1} g_2^{k_2} \cdots   g_{i}^{k_{i}} g_{i+1}^1 g_{i+2}^{k_{i+2}} \cdots g_n^{k_n} \bigr) \Bigl)
\end{array} \nonumber \\
&= \displaystyle \sum_{\bm{k} \in \mathcal{K}_{i+1}^1} (-1)^{|\bm{k}|-1}A_i = -\displaystyle \sum_{\bm{k} \in \mathcal{K}_{i+1}^1} (-1)^{ |\bm{k}| } A_i. \label{equation03}
\end{align}
}
Using \eqref{equation02} and \eqref{equation03}, we have
\begin{align}
\psi_{n-1} \big((\bm{g}_{\{i\}};\bm{m}_{\{i\}})- ( \bm{g}_{\{ \triangleleft i \}} ; \bm{m}_{\{ \triangleleft i \}})\big) &= \displaystyle \sum_{\bm{k} \in \mathcal{K}_{i+1}^0} (-1)^{ |\bm{k}| } A_i +\displaystyle \sum_{\bm{k} \in \mathcal{K}_{i+1}^1} (-1)^{ |\bm{k}| } A_i \nonumber \\
&= \displaystyle \sum_{\bm{k} \in \mathcal{K}_n} (-1)^{ |\bm{k}| } A_i. \label{equation04}
\end{align}
Using \eqref{cad}, \eqref{equality02} and \eqref{equation04}, we have
\begin{align*}
\psi_{n-1}\circ \partial_n^{\rm BR_U} (\bm{g}; \bm{m}) &= \displaystyle \sum_{i=1}^{n-1} (-1)^i \displaystyle \sum_{\bm{k} \in \mathcal{K}_n} (-1)^{|\bm{k}|} A_i = L.
\end{align*}

Secondly, we show $\partial_{n}^{\rm gp}\circ\psi_{n}(\bm{g};\bm{m})=L$.
We have
\begin{align}
\partial_{n}^{\rm gp}\circ\psi_{n}(\bm{g};\bm{m})
&= \displaystyle \sum_{\bm{k}\in\mathcal{K}_n}(-1)^{|\bm{k}|} \partial_{n}^{\rm gp}(m_1\bm{g}_{\bm{k}_1}, \ldots, m_n\bm{g}_{\bm{k}_n})\nonumber\\
&= \displaystyle \sum_{\bm{k}\in\mathcal{K}_n}(-1)^{|\bm{k}|}  (m_2\bm{g}_{\bm{k}_2}, \ldots, m_n\bm{g}_{\bm{k}_n})+L \nonumber \\
&\phantom{=} +\displaystyle \sum_{\bm{k}\in\mathcal{K}_n}(-1)^{|\bm{k}|}(-1)^n (m_1\bm{g}_{\bm{k}_1}, \ldots, m_{n-1}\bm{g}_{\bm{k}_{n-1}}). \label{aast}
\end{align}
We have
\begin{align}
&\displaystyle \sum_{\bm{k}\in\mathcal{K}_n}(-1)^{|\bm{k}|}  (m_2\bm{g}_{\bm{k}_2}, \ldots, m_n\bm{g}_{\bm{k}_n}) \nonumber \\
&= \displaystyle \sum_{\bm{k}\in\mathcal{K}^0_2}(-1)^{|\bm{k}|} (m_2\bm{g}_{\bm{k}_2}, \ldots, m_n\bm{g}_{\bm{k}_n}) + \displaystyle \sum_{\bm{k}\in\mathcal{K}^1_2}(-1)^{|\bm{k}|} (m_2\bm{g}_{\bm{k}_2}, \ldots, m_n\bm{g}_{\bm{k}_n}) \nonumber \\
&=0. \label{equation05}
\end{align}
This is because it holds that
\begin{align*}
&\displaystyle \sum_{\bm{k}\in\mathcal{K}^1_2}(-1)^{|\bm{k}|} \left( m_2\bm{g}_{\bm{k}_2},\ldots,m_n\bm{g}_{\bm{k}_n} \right) \\
&\begin{array}{l@{}l@{}l@{}}
= & \phantom{-} \displaystyle \sum_{\bm{k}\in\mathcal{K}^1_2}(-1)^{|\bm{k}|} & \left( m_2\varphi( g_1^{k_1} g_2^1 ) g_3^{k_3}\cdots g_n^{k_n},\ldots,m_n \varphi(g_1^{k_1} g_2^1g_3^{k_3}\cdots g_n^{k_n}) \right) \\
= & \phantom{-} \displaystyle \sum_{\bm{k}\in\mathcal{K}^1_2}(-1)^{|\bm{k}|} & \left( m_2 \varphi( g_1^{k_1} ) g_3^{k_3}\cdots g_n^{k_n},\ldots,m_n \varphi( g_1^{k_1} g_3^{k_3}\cdots g_n^{k_n}) \right) \\
= & \phantom{-} \displaystyle \sum_{\bm{k}\in\mathcal{K}^0_2}(-1)^{|\bm{k}|+1} & \left( m_2 \varphi( g_1^{k_1} ) g_2^0 g_3^{k_3} \cdots g_n^{k_n},\ldots,m_n \varphi(g_1^{k_1} g_2^0 g_3^{k_3}\cdots g_n^{k_n}) \right) \\
= & -\displaystyle \sum_{\bm{k}\in\mathcal{K}^0_2}(-1)^{|\bm{k}|} & \left( m_2\bm{g}_{\bm{k}_2},\ldots,m_n\bm{g}_{\bm{k}_n} \right).
\end{array}
\end{align*}
Similarly, we have
\begin{align}
&\displaystyle \sum_{\bm{k}\in\mathcal{K}_n}(-1)^{|\bm{k}|}(-1)^n (m_1\bm{g}_{\bm{k}_1}, \ldots, m_{n-1}\bm{g}_{\bm{k}_{n-1}}) \nonumber \\
&= (-1)^n \Bigl( \displaystyle \sum_{\bm{k}\in\mathcal{K}^0_n}(-1)^{|\bm{k}|} (m_1\bm{g}_{\bm{k}_1}, \ldots, m_{n-1} \bm{g}_{\bm{k}_{n-1}}) + \displaystyle \sum_{\bm{k}\in\mathcal{K}^1_n}(-1)^{|\bm{k}|} (m_1\bm{g}_{\bm{k}_1}, \ldots, m_{n-1} \bm{g}_{\bm{k}_{n-1}}) \Bigr) \nonumber \\
&=0. \label{equation06}
\end{align}
By  \eqref{aast}, \eqref{equation05} and \eqref{equation06}, we have $\partial_{n}^{\rm gp}\circ\psi_{n}(\bm{g};\bm{m})=L$.
\end{proof}
\subsection{Birack cocycles of $G$-Alexander biquandles}\label{subsection:A2}
Let $\gamma = (\gamma_n)$ and $\psi = (\psi_n)$ be the chain maps defined in 
Subsection \ref{subsection:A1}.  We have a sequence 
\[ C_n^{\rm BR} (X; \mathbb Z) \overset{\gamma_n}{\longrightarrow} C_{n}^{\rm BR_U} (X; \mathbb Z) \overset{\psi_n}{ \longrightarrow } C_{n}^{\rm gp} (M; \mathbb Z)_G \]
of chain groups $C_n^{\rm BR} (X; \mathbb Z), C_{n}^{\rm BR_U} (X; \mathbb Z),  C_{n}^{\rm gp} (M; \mathbb Z)_G$ and chain maps $\gamma, \psi$ for $n\geq 1$.

\begin{theorem}\label{birack_cocycle_main}
For any $n$-cocycle $f: C_{n}^{\rm gp} (M; \mathbb Z)_G \to A$ of $M$, the map \[ \Phi_f:=f \circ  \psi_n \circ \gamma_n  : C_n^{\rm BR} (X; \mathbb Z)  \to A \] is a birack $n$-cocycle of the $G$-Alexander biquandle $X$.
\end{theorem}

\begin{proof}
Since $\psi_n \circ \gamma_n$ is a chain map, we see the result. 
\end{proof}

An $A$-multilinear map $f:M^n \to A$ is {\it $G$-invariant } if
$f(m_1g, \ldots, m_ng) = f(m_1, \ldots, m_n)$ for any $g \in G$ and $(m_1, \ldots, m_n) \in M^n$. 
Note that any $G$-invariant $A$-multilinear map $f:M^n \to A$ induces a cocycle $f:C_{n}^{\rm gp} (M; \mathbb Z)_G \to A$.

The following theorem follows from a direct calculation.

\begin{theorem}\label{birack cocycle explicit}
\begin{itemize}
\item[(1)] 
Let $f: M^2 \to A$ be a $G$-invariant $A$-multilinear map. 
The birack $2$-cocycle $\Phi_f= f \circ  \psi_2 \circ \gamma_2   : C_2^{\rm BR} (X; \mathbb Z)  \to A$ is formulated as 
\[
\Phi_f((m_1,g_1), (m_2,g_2))=f\bigl(m_1-m_2,m_2(1-\varphi(g_2)g_2^{-1})\bigr)
\]
for $((m_1,g_1), (m_2,g_2)) \in X^2 \subset C_2^{\rm BR} (X; \mathbb Z)$.
\item[(2)] 
Let $f: M^3 \to A$ be a $G$-invariant $A$-multilinear map. 
The birack $3$-cocycle $\Phi_f= f \circ  \psi_3 \circ \gamma_3   : C_3^{\rm BR} (X; \mathbb Z)  \to A$ is
 formulated as 
\begin{align*}
&\Phi_f((m_1,g_1), (m_2,g_2), (m_3,g_3))\\
&=f\bigl((m_1-m_2)(1-\varphi(g_2)^{-1}g_2),m_2-m_3,m_3(1-\varphi(g_3)g_3^{-1})\bigr)
\end{align*}
for $((m_1,g_1), (m_2,g_2), (m_3,g_3)) \in X^3\subset C_3^{\rm BR} (X; \mathbb Z)$.
\end{itemize}
\end{theorem}
\section{Cocycles of $G$-Alexander biquandles with the $X$-set $X$}\label{app:12}

Throughout this section, let $X= M \times G$ be the $G$-Alexander biquandle of $(M, \varphi)$ and we assume that $X$ is also an $X$-set with the action $*:=\uline{*}$, that is, it holds that $(x_0*x_1)*(x_2 \oline{*} x_1)=(x_0*x_2)*(x_1 \uline{*} x_2)$ for any $x_0, x_1, x_2 \in X$.
We discuss  birack cocycles of $X$ with the $X$-set $X$.   

Since we apply a similar argument as shown in Section~\ref{app:11} to this case, we summarize all the properties without proof. 

We define a chain complex  of $X = M \times G$ with the $X$-set $X$, denoted by $C^{\rm BR_U}_\ast (X;\mathbb Z)_X$, and define 
chain maps $\gamma$ and  $\psi$. It turns out that the chain complex $C^{\rm BR_U}_\ast (X;\mathbb Z)_X$ is isomorphic to the birack chain complex 
$C^{\rm BR}_\ast (X;\mathbb Z)_X$ and hence the homology group 
$H^{\rm BR_U}_n (X;\mathbb Z)_X$ is isomorphic to $H^{\rm BR}_n (X;\mathbb Z)_X$.  

Our goal in this section is to give Theorems~\ref{birack_cocycle_with_X-set_main} and~\ref{birack_cocycle_explicit_Xset}.

\subsection{The chain complex $C^{\rm BR_U}_\ast (X;\mathbb Z)_X$ and chain maps $\gamma$ and $\psi$}\label{subsection:B1}

\subsubsection{The chain complex $C^{\rm BR_U}_\ast (X;\mathbb Z)_X$ and the 
homology group $H^{\rm BR_U}_n(X; \mathbb Z)_X$}

For $\bm{g}= (g_0, g_1,  \ldots , g_n) \in G \times G^n$ and $\bm{m} = (m_0, m_1, \ldots, m_n) \in M \times M^n$, 
we use the following notations:
\begin{align*}
\bm{g}_{\{i\}} &:= (g_0, g_1,  \ldots, g_{i},g_{i+2},\ldots , g_n) \in G \times G^{n-1},\\
\bm{g}_{\{ \triangleleft i\} } &:= (g_{i+1}^{-1} g_0 g_{i+1},g_{i+1}^{-1} g_1 g_{i+1},  \ldots, g_{i+1}^{-1} g_{i} g_{i+1},g_{i+2},\ldots , g_n) \in G \times G^{n-1},\\
\bm{m}_{\{i\}} &:= (m_0, m_1, \ldots, m_{i-1},m_{i} + m_{i+1},m_{i+2}, \ldots , m_n) \in M \times M^{n-1} \text{ and}\\
\bm{m}_{\{ \triangleleft i\} } &:= \bigl( m_0g_{i+1}, m_1 g_{i+1}, \ldots, m_{i-1} g_{i+1},m_{i} g_{i+1} + m_{i+1} \varphi(g_{i+1}),\\
&\hspace{21.5ex}m_{i+2} \varphi(g_{i+1}), \ldots , m_n \varphi(g_{i+1}) \bigr) \in M \times M^{n-1},
\end{align*}
where $n$ and $i$ are integers with $n \ge 2$ and $0 \le i \le n-1$.

Let $C^{\rm BR_U}_n(X; \mathbb Z)_X$ be the free abelian group generated by the elements
\[
(\bm{g}, \bm{m})=(g_0, g_1, \ldots, g_n; m_0, m_1, \ldots, m_n)\in (G \times G^{n}) \times (M \times M^{n})
\]
 if $n\geq 1$, and $C^{\rm BR_U}_n(X; \mathbb Z)_X:=0$ otherwise.

We define a boundary map 
$\partial_n^{\rm BR_U}: C^{\rm BR_U}_n(X; \mathbb Z)_X \to C^{\rm BR_U}_{n-1}(X; \mathbb Z)_X$ by 
\[
\partial_n^{\rm BR_U} (\bm{g};\bm{m})
=\displaystyle \sum_{i=0}^{n-1} (-1)^i  \big\{(\bm{g}_{\{i\}};\bm{m}_{\{i\}})- ( \bm{g}_{\{ \triangleleft i \}};\bm{m}_{\{ \triangleleft i \}})\big\}
\]
if $n\geq 2$, and $\partial_n^{\rm BR_U}=0$ otherwise.
For example, we have
{\normalsize
\begin{align*}
&\partial_3^{\rm BR_U} (g_0, \, g_1, \, g_2, \, g_3; \, m_0, \, m_1, \, m_2, \, m_3)\\
&= (-1)^0 \bigl\{ (g_0,\, g_2, \, g_3; \, m_0+m_1, m_2, m_3)\\
&\phantom{=\bigl\{} -(g_1^{-1}g_0g_1, \, g_2, \, g_3; \, m_0g_1+m_1\varphi(g_1), \, m_2\varphi(g_1), \, m_3\varphi(g_1) )\bigr\}\\
&\phantom{=} +(-1)^1\bigl\{(g_0, \, g_1, \, g_3; \, m_0, m_1+m_2, m_3)\\
&\phantom{=\bigl\{} -(g_2^{-1}g_0g_2, \, g_2^{-1}g_1g_2, \, g_3; \, m_0g_2, \, m_1g_2+m_2\varphi(g_2), \, m_3\varphi(g_2)) \bigr\}\\
&\phantom{=} +(-1)^2 \bigl\{ (g_0, \, g_1, \, g_2; \, m_0, m_1, m_2+m_3)\\
&\phantom{=\bigl\{} -(g_3^{-1}g_0g_3, \, g_3^{-1}g_1g_3, \, g_3^{-1}g_2g_3; \, m_0g_3, \, m_1g_3, \, m_2g_3+m_3\varphi(g_3)) \bigr\}.
\end{align*}
}
\begin{lemma}
$C^{\rm BR_U}_*(X; \mathbb Z)_X=(C^{\rm BR_U}_n(X; \mathbb Z)_X, \partial_n^{\rm BR_U})_{n \in \mathbb Z}$ is a chain complex. 
\end{lemma}

\subsubsection{The chain map $\gamma: C_*^{\rm BR} (X; \mathbb Z)_X \to C_*^{\rm BR_U} (X; \mathbb Z)_X$}

For $n\geq 1$,  define a homomorphism $\gamma_n: C_n^{\rm BR} (X; \mathbb Z)_X \to C_{n}^{\rm BR_U} (X; \mathbb Z)_X$
by 
$$
\gamma_n((m_0, g_0), (m_1, g_1), \ldots , (m_n, g_n) ) := (\bm{g}; m_0', m_1', \ldots, m_{n-1}', m_n),
$$
where we write $m_i':=m_i-m_{i+1}$ and $\bm{g} = (g_0, \dots, g_n)$.
Define $\gamma_n=0$ for  $i<1$.

\begin{lemma}\label{lem-gamchain}
\begin{itemize}
\item[(1)]
The map $\gamma$ is a chain map, that is, it holds that
\[
\gamma_{n-1} \circ \partial_n^{\rm BR}= \partial_n^{\rm BR_U} \circ \gamma_{n}.
\]
\item[(2)] For each integer $n$, the map $\gamma_n$ is an isomorphism. 
\end{itemize} 
\end{lemma}

\subsubsection{The chain map $\psi : C_{*}^{\rm BR_U} (X; \mathbb Z)_X \to C_{*+1}^{\rm gp} (M; \mathbb Z)_G$}

For $n\geq 1$, define $\psi_n : C_{n}^{\rm BR_U} (X; \mathbb Z)_X \to C_{n+1}^{\rm gp} (M; \mathbb Z)_G$ by 
\[
\begin{array}{l}
\psi_n (g_0, g_1, \ldots, g_n; m_0, m_1, \ldots, m_n) := \displaystyle \sum_{\bm{k} \in \mathcal K_n} (-1)^{|\bm{k}|} (m_0\bm{g}_{\bm{k}_0},m_1 \bm{g}_{\bm{k}_1}, \ldots, m_n \bm{g}_{\bm{k}_n}),
\end{array}
\]
where $\mathcal{K}_n :=\{0,1\}^n$ and 
\begin{align*}
         |\bm{k}| &:= k_1 + k_2 + \cdots  +k_n,\\
 \bm{g}_{\bm{k}_0}&:= g_1^{k_1}g_2^{k_2}\cdots g_n^{k_n} \in G \text{ and}\\
\bm{g}_{\bm{k}_i} &:= \varphi(g_1^{k_1} g_2^{k_2} \cdots g_{i}^{k_i}) g_{i+1}^{k_{i+1}}g_{i+2}^{k_{i+2}} \cdots g_{n}^{k_{n}} \in G
\end{align*}
for an element $\bm{k}=(k_1, \ldots, k_n) \in \mathcal K_n$ and an integer $i$ with $1 \le i \le n$.
Define $\psi_n=0$ for $n <1$.
We note that the codomain of $\psi_n$ is $C_{n+1}^{\rm gp} (M; \mathbb Z)_G$.

For example, we have
\begin{align*}
\psi_3(\bm{g};\bm{m}) &=(-1)^{0+1+1}(m_0g_1^0g_2^1g_3^1, m_1\varphi(g_1^0)g_2^1g_3^1, m_2\varphi(g_1^0g_2^1)g_3^1, m_3\varphi(g_1^0g_2^1g_3^1) )\\
&\phantom{=}+(-1)^{0+1+0}(m_0g_1^0g_2^1g_3^0, m_1\varphi(g_1^0)g_2^1g_3^0, m_2\varphi(g_1^0g_2^1)g_3^0, m_3\varphi(g_1^0g_2^1g_3^0) )\\
&\phantom{=}+(-1)^{0+0+1}(m_0g_1^0g_2^0g_3^1, m_1\varphi(g_1^0)g_2^0g_3^1, m_2\varphi(g_1^0g_2^0)g_3^1, m_3\varphi(g_1^0g_2^0g_3^1) )\\
&\phantom{=}+(-1)^{0+0+0}(m_0g_1^0g_2^0g_3^0, m_1\varphi(g_1^0)g_2^0g_3^0, m_2\varphi(g_1^0g_2^0)g_3^0, m_3\varphi(g_1^0g_2^0g_3^0) )\\
&\phantom{=}+(-1)^{1+1+1}(m_0g_1^1g_2^1g_3^1, m_1\varphi(g_1^1)g_2^1g_3^1, m_2\varphi(g_1^1g_2^1)g_3^1, m_3\varphi(g_1^1g_2^1g_3^1) )\\
&\phantom{=}+(-1)^{1+1+0}(m_0g_1^1g_2^1g_3^0, m_1\varphi(g_1^1)g_2^1g_3^0, m_2\varphi(g_1^1g_2^1)g_3^0, m_3\varphi(g_1^1g_2^1g_3^0) )\\
&\phantom{=}+(-1)^{1+0+1}(m_0g_1^1g_2^0g_3^1, m_1\varphi(g_1^1)g_2^0g_3^1, m_2\varphi(g_1^1g_2^0)g_3^1, m_3\varphi(g_1^1g_2^0g_3^1) )\\
&\phantom{=}+(-1)^{1+0+0}(m_0g_1^1g_2^0g_3^0, m_1\varphi(g_1^1)g_2^0g_3^0, m_2\varphi(g_1^1g_2^0)g_3^0, m_3\varphi(g_1^1g_2^0g_3^0) ).
\end{align*}

\begin{lemma}
The map $\psi = (\psi_n)$ is a chain map, that is, it holds that
\[
\psi_{n-1} \circ \partial_n^{\rm BR_U}= \partial_{n+1}^{\rm gp} \circ \psi_{n}.
\]
\end{lemma}
\subsection{Cocycles of $G$-Alexander biquandles with the $X$-set $X$}
As a consequence of Subsection \ref{subsection:B1}, we have a sequence 
\[
C_n^{\rm BR} (X; \mathbb Z)_X \overset{\gamma_n}{\longrightarrow} C_{n}^{\rm BR_U} (X; \mathbb Z)_X \overset{\psi_n}{ \longrightarrow } C_{n+1}^{\rm gp} (M; \mathbb Z)_G
\]
of chain groups $C_n^{\rm BR} (X; \mathbb Z)_X, C_{n}^{\rm BR_U} (X; \mathbb Z)_X,  C_{n}^{\rm gp} (M; \mathbb Z)_G$ and chain maps $\gamma, \psi$ for $n\geq 1$.
Therefore we have the following theorem.

\begin{theorem}\label{birack_cocycle_with_X-set_main}
For any $(n+1)$-cocycle $f: C_{n+1}^{\rm gp} (M; \mathbb Z)_G \to A$, the map
$$\Phi_f:=f \circ  \psi_n \circ \gamma_n  : C_n^{\rm BR} (X; \mathbb Z)_X  \to A$$  
is a birack $n$-cocycle of the $G$-Alexander biquandle $X=M \times G$.
\end{theorem}

\begin{theorem}\label{birack_cocycle_explicit_Xset}
\begin{itemize}
\item[(1)] 
Let $f: M^3 \to A$ be a $G$-invariant $A$-multilinear map. 
The birack $2$-cocycle $\Phi_f= f \circ  \psi_2 \circ \gamma_2   : C_2^{\rm BR} (X; \mathbb Z)_X  \to A$ of the $G$-Alexander biquandle $X=M \times G$ is formulated as 
\begin{align*}
&\Phi_f((m_0,g_0),(m_1,g_1), (m_2,g_2))\\
&=f \bigl( m_0' (1-\varphi(g_1)^{-1}g_1),m_1', m_2(1-\varphi(g_2)g_2^{-1}) \bigr)
\end{align*}
for $((m_0,g_0),(m_1,g_1), (m_2,g_2)) \in X\times X^2 \subset C_2^{\rm BR} (X; \mathbb Z)_X$,
where $m_i':=m_i -m_{i+1}$.
\item[(2)] 
Let $f: M^4 \to A$ be a $G$-invariant $A$-multilinear map. 
The birack $3$-cocycle $\Phi_f= f \circ  \psi_3 \circ \gamma_3   : C_3^{\rm BR} (X; \mathbb Z)_X  \to A$ of the $G$-Alexander biquandle $X=M \times G$ is
 formulated as 
\begin{align*}
&
\hspace{-2em}
\Phi_f \bigl( (m_0,g_0), (m_1,g_1), (m_2,g_2), (m_3,g_3) \bigl)\\
&
\hspace{-1.5em}
=f \bigl( m_0' (1-\varphi(g_1)^{-1}g_1), m_1', m_2', m_3(1-\varphi(g_3)g_3^{-1}) \bigr)  \\
&
\hspace{-1.5em}
\phantom{=} - f \bigl( m_0' (1-\varphi(g_1)^{-1}g_1)g_2, m_1' g_2, m_2' \varphi(g_2), m_3(1-\varphi(g_3)g_3^{-1})\varphi(g_2) \bigr)
\end{align*}
for $((m_0,g_0), (m_1,g_1), (m_2,g_2), (m_3,g_3)) \in X \times X^3 \subset C_3^{\rm BR} (X; \mathbb Z)_X$,
where $m_i':=m_i -m_{i+1}$.
\end{itemize}
\end{theorem}
\section{Cocycles of $G$-Alexander multiple conjugation biquandles}\label{app:21}

Throughout this section, let $X=\bigsqcup_{m \in M}(\{m\} \times G)=M \times G$ be the $G$-Alexander multiple conjugation biquandle of $(M, \varphi)$, see Definition \ref{def_Gfamily_of_Alexander_quandle}.
Our goal in this section is to give Theorem~\ref{mcb_cocycles_no_X-set_main}. 

\subsection{Degenerate subcomplexes
$D^{\rm BR}_*(X; \mathbb Z)$, $D^{\rm BR_U}_*(X;\mathbb Z)$ and the induced homomorphisms $\gamma_n$, $\psi_{n,\lambda}$}\label{subsection:A3}

\subsubsection{The degenerate subcomplex $D^{\rm BR}_*(X;\mathbb Z)$ of $C^{\rm BR}_*(X;\mathbb Z)$}

Let $D^{\rm BR}_n(X;\mathbb Z)$ be the subgroup of $C^{\rm BR}_n(X;\mathbb Z)$ generated by the elements of the following two sets
\begin{align*}
&
\displaystyle \bigcup_{ i=1}^{ n-1} \Big\{ ( \bm{x}^{i-1} , (m, g),  (m, h),  \bm{x}_{i+2}) ~\Big|~\bm{x} \in X^n,~m \in M,~g,h \in G \Big\} \text{ and} \\
&
\displaystyle \bigcup_{i=1}^{n} \left\{
\begin{array}{lcl}
\begin{minipage}{6.6cm}
{\normalsize
$\hspace{1.8ex}(\bm{x}^{i-1} , (m, gh),  \bm{x}_{i+1} ) - ( \bm{x}^{i-1} , (m, g),  \bm{x}_{i+1})$\\
$ - \bigl( \bm{x}^{i-1} \underline{*} (m,g),~ \bigl( (m, h),  \bm{x}_{i+1} \bigr) \overline{*} (m,g) \bigl)$}
\end{minipage}
& \Bigg| &
\begin{minipage}{2.7cm}
{\normalsize
$\bm{x} \in X^n,\\ m \in M,~g,h \in G$
}
\end{minipage}
\end{array}
\right\}
\end{align*}
for $n \geq 2$.
We define $D^{\rm BR}_n(X;\mathbb Z)=0$ for $n \leq 1$.
We note that 
\begin{align*}
&\left( \bm{x}^{i-1} \underline{*} (m,g),\bigl( (m, h),  \bm{x}_{i+1} \bigr) \overline{*} (m,g) \right)\\
&=(x_1\uline*(m,g), \ldots, x_{i-1}\uline*(m,g), (m,h)\oline*(m,g),x_{i+1}\oline*(m,g), \ldots, x_n\oline*(m,g)).
\end{align*}
\begin{lemma} 
$D^{\rm BR}_*(X;\mathbb Z):=(D^{\rm BR}_n(X;\mathbb Z), \partial_n^{\rm BR})_{n \in \mathbb Z}$ is a subcomplex of $C^{\rm BR}_*(X;\mathbb Z)$.
\end{lemma}
\begin{proof}
We fix an integer $n \ge 2$ and show $\partial_n^{\rm BR} ( D^{\rm BR}_n(X;\mathbb Z) ) \subset D^{\rm BR}_{n-1}(X;\mathbb Z)$.
It suffices to show
\[ \partial_n^{\rm BR}(\bm{x}^{i-1},(m,g),(m,h),\bm{x}_{i+2})\equiv 0 \] for any $i$ with $1 \le i \le n-1$ and 
\begin{align*}
&\partial_n^{\rm BR}(\bm{x}^{i-1},(m,gh),\bm{x}_{i+1}) \equiv \partial_n^{\rm BR}(\bm{x}^{i-1},(m,g),\bm{x}_{i+1})\\
&\hspace{28ex}+ \partial_n^{\rm BR} \bigl( \bm{x}^{i-1} \underline{*}(m,g), \bigl( (m,h),\bm{x}_{i+1} \bigr) \overline{*}(m,g) \bigr)
\end{align*}
for any $i$ with $1 \le i \le n$ in $C_{n-1}^{\rm BR}(X;\Bbb Z)/D_{n-1}^{\rm BR}(X;\Bbb Z)$.
\vspace{1ex}

We verify the first equality in the quotient group. Put $x_i:=(m,g)$, $x_{i+1}:=(m,h)$ and $\bm{x}:=(x_1, \ldots, x_n)$.
\begin{align}
&\partial_n^{\rm BR}(\bm{x}) \nonumber \\
&=\displaystyle \sum_{j=1}^{n} (-1)^{j-1} \bigl\{ (\bm{x}^{j-1}, \bm{x}_{j+1})-( \bm{x}^{j-1} \uline{*} x_j , \bm{x}_{j+1} \oline{*} x_j) \bigr\} \nonumber \\
&= \displaystyle \sum_{j=1}^{i-1} (-1)^{j-1} \underline{\{ (\bm{x}^{j-1}, \bm{x}_{j+1}) -(\bm{x}^{j-1} \uline{*}x_j , \bm{x}_{j+1}\oline{*} x_j)\}}_{(A)} \nonumber \\
&\phantom{=} + (-1)^{i-1} \underline{\bigl\{ (\bm{x}^{i-1},(m,h),\bm{x}_{i+2})- \bigl( \bm{x}^{i-1}\underline{*}(m,g), ( (m,h),\bm{x}_{i+2} ) \oline{*}(m,g) \bigr) \bigr\}}_{(B)} \nonumber  \\
&\phantom{=} + (-1)^{i \ \ \ } \underline{\bigl\{ (\bm{x}^{i-1},(m,g),\bm{x}_{i+2})- \bigl( ( \bm{x}^{i-1}, (m,g) ) \uline{*}(m,h), \bm{x}_{i+2} \overline{*}(m,h) \bigr) \bigr\}}_{(C)} \nonumber \\
&\phantom{=} + \displaystyle \sum_{j=i+2}^{n} (-1)^{j-1} \underline{ \{ (\bm{x}^{j-1}, \bm{x}_{j+1}) -(\bm{x}^{j-1} \uline{*}x_j , \bm{x}_{j+1}\oline{*} x_j)\}}_{(D)}. \label{ewq1}
\end{align}

Since $(m,gh)=(m, h(h^{-1}gh))$ and $(m,g)\uline{*}(m,h) = (mh+m(\varphi(h)-h), h^{-1}gh) = (m, h^{-1}gh)\oline{*}(m,h)$, we see that 
for $(B)$ and $(C)$ of \eqref{ewq1}, 
\begin{align*}
&(B)-(C) \\
&=(\bm{x}^{i-1},(m,h),\bm{x}_{i+2})-(\bm{x}^{i-1} \underline{*}(m,g),( (m,h) , \bm{x}_{i+2})  \oline{*} (m,g) ) -(\bm{x}^{i-1},(m,g),\bm{x}_{i+2})\\
&\phantom{=} + \bigl( ( \bm{x}^{i-1}, (m,g) ) \uline{*} (m,h), \bm{x}_{i+2} \overline{*} (m,h) \bigr) \\
&\equiv  (\bm{x}^{i-1},(m,h),\bm{x}_{i+2}) -(\bm{x}^{i-1},(m,gh),\bm{x}_{i+2})\\
&\phantom{=} + \bigl( ( \bm{x}^{i-1}, (m,g) ) \uline{*} (m,h), \bm{x}_{i+2} \overline{*} (m,h) \bigr) \\
&= (\bm{x}^{i-1},(m,h),\bm{x}_{i+2}) -(\bm{x}^{i-1},(m, h(h^{-1}gh)),\bm{x}_{i+2})\\
&\phantom{=} + (\bm{x}^{i-1}\uline{*}(m,h), ( (m, h^{-1}gh),\bm{x}_{i+2} ) \oline{*} (m,h) )\\
&\equiv 0.
\end{align*}


When $j\not=i$ and $j\not=i+1$, we have $(\bm{x}_{j-1}, \bm{x}^{j+1}) \equiv 0$ immediately. Since the first element of $(m,g)\uline* x_j$ is equal to that of $(m,h)\uline* x_j$ and the first element of $(m,g)\oline* x_j$ is equal to that of $(m,h)\oline* x_j$, we have
\begin{align*}
\left( \bm{x}^{j-1}\underline{*}x_j, \left( \bm{x}_{j+1}^{i-1}, (m,g), (m,h), \bm{x}^{i+2} \right) \overline{*}x_j \right)
&\in D_{n-1}^{\rm BR}(X;\Bbb Z) \text{ and}\\
\left( \bigl( \bm{x}^{i-1},(m,g),(m,h),\bm{x}_{i+2}^{j-1} \bigr) \underline{*} x_j, \bm{x}^{j+1} \overline{*}x_j \right)
&\in D_{n-1}^{\rm BR}(X;\Bbb Z),
\end{align*}
where $\bm{x}_{a}^{b}$ means 
the sequence $x_{a}, x_{a+1}, \ldots, x_{b-1}, x_{b}$. Then, $(\bm{x}_{j-1} \uline{*}x_j , \bm{x}^{j+1}\oline{*} x_i)\equiv0$ in the equality of \eqref{ewq1}. Then, we have $(A)\equiv (D)\equiv 0$.
Hence, we have $\partial_n^{\rm BR}(\bm{x})\equiv0$.
\vspace{1ex}

We verify the second equality in the quotient group. Put $x_i:=(m,gh)$ and $\bm{x}:=(x_1, \ldots, x_n)=(x_1, \ldots, x_{i-1} ,(m,gh), x_{i+1}, \ldots, x_n)$. We have
\begin{align*}
\partial_n^{\rm BR}(\bm{x}) &= \displaystyle \sum_{j=1}^{n} (-1)^{j-1} \bigl\{ (\bm{x}^{j-1}, \bm{x}_{j+1})-(\bm{x}^{j-1} \uline{*}x_j , \bm{x}_{j+1}\oline{*} x_j) \bigr\}\\
&= \phantom{+} \displaystyle \sum_{j=1}^{i-1}(-1)^{j-1} \bigl\{ (\bm{x}^{j-1}, \bm{x}_{j+1})-(\bm{x}^{j-1} \uline{*}x_j , \bm{x}_{j+1}\oline{*} x_j)\bigr\}\\
&\phantom{=} + \displaystyle \sum_{j=i}^{i} (-1)^{j-1} \bigl\{ (\bm{x}^{j-1}, \bm{x}_{j+1})-(\bm{x}^{j-1} \uline{*}x_j , \bm{x}_{j+1}\oline{*} x_j)\bigr\}\\
&\phantom{=} + \displaystyle \sum_{j=i+1}^{n} (-1)^{j-1} \bigl\{ (\bm{x}^{j-1}, \bm{x}_{j+1})-(\bm{x}^{j-1} \uline{*}x_j , \bm{x}_{j+1}\oline{*} x_j)\bigr\}.
\end{align*}

\noindent
When $1\le j \le i-1$, we have
\begin{align*}
( \bm{x}^{j-1}, \bm{x}_{j+1} ) &= ( \bm{x}^{j-1},\bm{x}_{j+1}^{i-1},(m,gh),\bm{x}_{i+1}) \\
&\equiv \underline{\bigl( \bm{x}^{j-1},\bm{x}_{j+1}^{i-1}, (m,g),\bm{x}_{i+1} \bigr)}_{(A)} \\
&\phantom{\equiv} + \underline{ \left(\left(\bm{x}^{j-1},\bm{x}_{j+1}^{i-1}\right) \uline{*} (m,g), \left( (m,h),\bm{x}_{i+1} \right) \oline{*} (m,g) \right)}_{(B)}.
\end{align*}

\noindent
Put $x_j:=(m_j, g_j)$.
Since $(m,g')\oline{*}x_j=(m\varphi(g_j), g')$ for any $g' \in G$, we have
\begin{align*}
&(\bm{x}^{j-1} \uline{*} x_j , \bm{x}_{j+1}\oline{*} x_j)\nonumber \\
&=( \bm{x}^{j-1}\uline{*}x_j, \left( \bm{x}_{j+1}^{i-1}, (m,gh),\bm{x}_{i+1} \right) \oline{*}x_j )\nonumber \\
&=( \bm{x}^{j-1}\uline{*}x_j,\bm{x}_{j+1}^{i-1}\oline{*}x_j, (m\varphi(g_j), gh),\bm{x}_{i+1}\oline{*}x_j ) \nonumber \\
&\equiv \bigl( \bm{x}^{j-1} \uline{*}x_j, \bm{x}_{j+1}^{i-1}\oline{*}x_j, ( m \varphi(g_j), g ),\bm{x}_{i+1}\oline{*}x_j \bigr) \nonumber \\
&\phantom{=} + \left( \left( \bm{x}^{j-1}\uline{*}x_j, \bm{x}_{j+1}^{i-1}\oline{*}x_j \right) \uline{*} \left( m \varphi(g_j), g \right), \left( \left( m \varphi(g_j), h \right), \bm{x}_{i+1} \oline{*}x_j \right) \oline{*} \left( m \varphi(g_j), g \right) \right) \nonumber \\
&=\underline{( \bm{x}^{j-1}\uline{*}x_j, \left( \bm{x}_{j+1}^{i-1}, (m,g), \bm{x}_{i+1} \right) \oline{*}x_j )}_{(C)} \label{qqq3} \\
&\phantom{=} + \underline{ \left( \left( \bm{x}^{j-1}\uline{*}x_j, \bm{x}_{j+1}^{i-1}\oline{*}x_j \right) \uline{*} \left( (m,g)\oline{*}x_j \right), \left( (m,h)\oline{*}x_j, \bm{x}_{i+1} \oline{*}x_j \right) \oline{*} \left( (m,g)\oline{*}x_j \right) \right)}_{(D)}.
\end{align*}

Then, we have
\begin{align*}
&(D) =\underline{ \Bigl( \left( \bm{x}^{j-1} \uline{*} (m,g) \right) \uline* (x_j \uline* (m,g)),\left( \bm{x}_{j+1}^{i-1} \uline{*}(m,g), (m,h) \oline* (m,g), \bm{x}_{i+1}^{}\oline* (m,g) \right) \oline* (x_j \uline* (m,g)) \Bigr)}_{(D')}.
\end{align*}

\noindent
When $j=i$, we have
\begin{align*}
&( \bm{x}^{j-1}, \bm{x}_{j+1} ) = \underline{( \bm{x}^{i-1}, \bm{x}_{i+1} )}_{(E)},  \\
&(\bm{x}^{j-1} \uline{*} x_j , \bm{x}_{j+1}\oline{*} x_j) \\
& =(\bm{x}^{i-1} \uline{*} (m,gh) , \bm{x}_{i+1}\oline{*} (m,gh)) \\
&=  \underline{\left( \left( \bm{x}^{i-1} \uline{*} (m,g) \right) \uline{*} \left( (m,h)\oline{*}(m,g) \right), \left( \bm{x}_{i+1}\oline* ( m, g ) \right) \oline{*} \left( (m,h)\oline{*} (m, g) \right) \right)}_{(F)}.
\end{align*}

\noindent
When $i+1 \le j \le n$, we have

\begin{align}
( \bm{x}^{j-1}, \bm{x}_{j+1} ) &=( \bm{x}^{i-1}, (m,gh), \bm{x}_{i+1}^{j-1},\bm{x}_{j+1} ) \nonumber \\
&\equiv \underline{( \bm{x}^{i-1}, (m,g), \bm{x}_{i+1}^{j-1},\bm{x}_{j+1} )}_{(G)} \nonumber \\
&\phantom{\equiv} + \underline{( \bm{x}^{i-1} \uline{*}(m,g),~\{(m,h), \bm{x}_{i+1}^{j-1},\bm{x}_{j+1} \}\oline{*}(m,g) )}_{(H)}.
\end{align}

\noindent
Put $x_j:=(m_j, g_j)$. We have $(m', g') \uline* x_j =( m' g_j + m_j ( \varphi(g_j) - g_j ) , g_j^{-1} g' g_j ) $ and $(m', g') \oline* m_j =( m' \varphi(g_j)  , g' ) $ for any $(m', g') \in X$. 
Put 
$A:=mg_j+m_j(\varphi(g_j)-g_j)$.
We have
\begin{align}
&(\bm{x}^{j-1} \uline{*} x_j , \bm{x}_{j+1}\oline{*} x_j)\nonumber \\
&= ( \bm{x}^{i-1} \uline* x_j, (A, g_j^{-1} g g_j g_j^{-1} h g_j), \bm{x}_{i+1}^{j-1}\uline{*} x_j, \bm{x}_{j+1}\oline{*} x_j )\nonumber \\
&\equiv ( \bm{x}^{i-1} \uline* x_j, (A, g_j^{-1} g g_j ), \bm{x}_{i+1}^{j-1}\uline{*} x_j, \bm{x}_{j+1}\oline{*} x_j ) \nonumber \\
&\phantom{=}+\bigl( \bigl( \bm{x}^{i-1} \uline* x_j \bigr) \uline* ( A , g_j^{-1} g g_j ) , \Bigl( ( A , g_j^{-1} h g_j ), \bm{x}_{i+1}^{j-1}\uline{*} x_j, \bm{x}_{j+1}\oline{*} x_j \Bigr) \oline* ( A, g_j^{-1} g g_j) \bigr) \nonumber \\
&= \underline{\bigl( \bigl( \bm{x}^{i-1} , (m,g) , \bm{x}_{i+1}^{j-1} \bigr) \uline* x_j, \bm{x}_{j+1} \oline* x_j \bigr)}_{(I)}  \nonumber \\
&\phantom{=}+\underline{\bigl( ( \bm{x}^{i-1} \uline* x_j ) \uline* \bigl( (m,g) \uline* x_j \bigr), \Bigl( (m, h) \uline* x_j, \bm{x}_{i+1}^{j-1} \uline* x_j, \bm{x}_{j+1}\oline* x_j \Bigr) \oline* \bigl( (m,g) \uline* x_j \bigr) \bigr)}_{(J)}. \label{qqq10}
\end{align}

Then, we have
\begin{align*}
&J= \nonumber \\
&\underline{\Bigl( \bigl( \bm{x}^{i-1} \uline{*}(m,g), (m,h)\oline* (m,g), \bm{x}_{i+1}^{j-1}\oline* (m,g) \bigr) \uline* (x_j\oline* (m,g)),\bigl( \bm{x}_{j+1} \oline{*}(m,g) \bigr) \oline* (x_j\oline* (m,g)) \Bigr)}_{(J')}.
\end{align*}
Hence, we have
\begin{align*}
\partial_n^{\rm BR}&(\bm{x}) \equiv \displaystyle \sum_{j=1}^{i-1}(-1)^{j-1} ( A+B-( C+D' ) )\\
& +(-1)^{i-1}(E-F) +\displaystyle \sum_{j=i+1}^{n}(-1)^{j-1} ( G+H- \left( I + J' \right) ).
\end{align*}
By the definition of the map $\partial_n^{\rm BR}$, we have
\begin{align*}
\partial_n^{\rm BR}&(\bm{x}^{i-1},(m,g),\bm{x}_{i+1}) = \displaystyle \sum_{j=1}^{i-1}(-1)^{j-1}(A- C)\\
& +(-1)^{i-1} ( E-(\bm{x}^{i-1}\uline* (m,g), \bm{x}_{i+1} \oline* (m,g) ) ) +\displaystyle \sum_{j=i+1}^{n}(-1)^{j-1}(G-I)
\end{align*}
and
\begin{align*}
\partial_n^{\rm BR}& (\bm{x}^{i-1}\underline{*}(m,g), \left( (m,h),\bm{x}_{i+1} \right) \overline{*}(m,g)) = \displaystyle \sum_{j=1}^{i-1}(-1)^{j-1} (B-D' ) \\
& +(-1)^{i-1} ( ( \bm{x}^{i-1}\underline{*}(m,g),\bm{x}_{i+1}\overline{*}(m,g) ) -F ) +\displaystyle \sum_{j=i+1}^{n}(-1)^{j-1} (H-J').
\end{align*}
Therefore, it holds that
\begin{align*}
& \partial_n^{\rm BR}(\bm{x})\equiv \partial_n^{\rm BR}(\bm{x}^{i-1},(m,g),\bm{x}_{i+1})+ \partial_n^{\rm BR}(\bm{x}^{i-1}\underline{*}(m,g), ( (m,h),\bm{x}_{i+1} ) \overline{*}(m,g)). 
\end{align*}
This completes the proof.  
\end{proof}

The \textit{normalized birack chain complex} is  $C^{\rm norBR}_*(X;\mathbb Z):=C^{\rm BR}_*(X;\mathbb Z)/D^{\rm BR}_*(X;\mathbb Z)$.  It determines 
the \textit{normalized birack homology group} $H^{\rm nor BR}_n (X;\mathbb Z)$.
In the ordinary way, for an abelian group $A$, we have the (co)homology theory with the coefficient group $A$ and the homology group $H^{\rm nor BR}_n (X;A)$ and the cohomology group $H_{\rm nor BR}^n (X;A)$ are defined.

\subsubsection{The degenerate subcomplex $D^{\rm BR_U}_*(X;\mathbb Z)$ of $C^{\rm BR_U}_*(X;\mathbb Z)$} 

We introduce the degenerate subcomplex $D^{\rm BR_U}_*(X;\mathbb Z)$ of $C^{\rm BR_U}_*(X;\mathbb Z)$, which is a counterpart of the degenerate subcomplex 
$D^{\rm BR}_*(X;\mathbb Z)$ of $C^{\rm BR}_*(X;\mathbb Z)$.  

Let $D^{\rm BR_U}_n(X;\mathbb Z)$ be the subgroup of $C^{\rm BR_U}_n(X;\mathbb Z)$ generated by the elements of the following sets
\begin{align*}
\displaystyle \bigcup_{ i=1}^{ n-1} &\Big\{ (\bm{g};\bm{m}^{i-1},0,\bm{m}_{i+1}) ~\Big|~ \bm{g} \in G^n,~ \bm{m} \in M^n\Big\} \text{ and}\\
\displaystyle \bigcup_{i=1}^{n} &\Bigl\{
\begin{array}{lcl}
\begin{minipage}{6.9cm}
{\normalsize
$ \hspace{2ex}(\bm{g}^{i-1}, g_i h, \bm{g}_{i+1};\bm{m} )$\\
$- (\bm{g};\bm{m} )$$-(g_i^{-1} \bm{g}^{i-1}g_i, h, \bm{g}_{i+1} ;  \bm{m}^{i-1}g_i, \bm{m}_i \varphi(g_i))$}
\end{minipage}
& \Big| &
\begin{minipage}{2.4cm}
{\normalsize
$h \in G$,
$\bm{g} \in G^n$,\\
$\bm{m} \in M^n$
}
\end{minipage}
\end{array}
\Bigr\}
\end{align*}
for  $n \geq 2$, where we write
\begin{align*}
&(\bm{g}; \bm{m}^{i-1},0,\bm{m}_{i+1}) := (g_1, \ldots, g_n;m_1, \ldots, m_{i-1}, 0, m_{i+1}, \ldots, m_n ),\\
&(\bm{g}^{i-1}, g_ih, \bm{g}_{i+1} ; \bm{m}):= (g_1, \ldots, g_{i-1}, g_i h, g_{i+1}, \ldots , g_n ; m_1, \ldots, m_n ) \text{ and}\\
&(g_i^{-1} \bm{g}^{i-1} g_i, h, \bm{g}_{i+1} ; \bm{m}^{i-1} g_i , \bm{m}_i \varphi(g_i) ):=( ( g_{i}^{-1} g_1 g_i ), \ldots, ( g_i^{-1} g_{i-1} g_i ),\\
&\hspace{23ex}h, g_{i+1}, \ldots, g_n; m_1g_i, \ldots, m_{i-1}g_i, m_i\varphi(g_i), \ldots, m_n\varphi(g_i)).
\end{align*}
We define $D^{\rm BR_U}_n(X;\mathbb Z)=0$ for $n \le1$.

\begin{lemma}\label{3lem_nor}
$D^{\rm BR_U}_*(X;\mathbb Z):=(D^{\rm BR_U}_n(X;\mathbb Z), \partial_n^{\rm BR_U})_{n \in \mathbb Z}$ is a subcomplex of $C^{\rm BR_U}_*(X;\mathbb Z)$.
\end{lemma}

\begin{proof}
In Lemma~\ref{3lem-resgam}, we will show that the isomorphism $\gamma_n $ defined in Subsection~\ref{subsection:A1} gives an isomorphism $D_n^{\rm BR} (X; \mathbb Z) \cong D_{n}^{\rm BR_U} (X; \mathbb Z)$.
\end{proof}

The \textit{ normalized U-birack chain complex} is  $C^{\rm norBR_U}_*(X;\mathbb Z):=C^{\rm BR_U}_*(X;\mathbb Z)/D^{\rm BR_U}_*(X;\mathbb Z)$. 
It determines the homology group $H^{\rm nor BR_U}_n(X;\mathbb Z)$. 
In the ordinary way, for an abelian group $A$, $H^{\rm nor BR_U}_n (X;A)$ and $H_{\rm nor BR_U}^n (X;A)$ are defined. 

\subsubsection{The induced homomorphism $\gamma_n$}

Next lemma shows that the isomorphism 
$\gamma_n: C_n^{\rm BR} (X; \mathbb Z) \to C_{n}^{\rm BR_U} (X; \mathbb Z)$, which is defined in Subsection~\ref{subsection:A1}, induces the isomorphism $\gamma_n:C_n^{\rm norBR} (X; \mathbb Z) \to C_{n}^{\rm norBR_U} (X; \mathbb Z)$, where we denote it by the same symbol  $\gamma_n$ for simplicity.

\begin{lemma}\label{3lem-resgam} 
It holds that $\gamma_n(D_{n}^{\rm BR} (X; \mathbb Z)) = D_{n}^{\rm BR_U} (X; \mathbb Z)$. 
Therefore,
$\gamma_n$ induces the isomorphism
$$\gamma_n: C_n^{\rm norBR} (X; \mathbb Z) \to C_{n}^{\rm norBR_U} (X; \mathbb Z).$$
\end{lemma} 

\begin{proof}
Let $\bm{x} \in  D_{n}^{\rm BR} (X; \mathbb Z)$.  If 
$\bm{x}=( \bm{x}^{i-1} , (m, g_i),  (m, g_{i+1}),  \bm{x}_{i+2})$,
 then we put $\gamma_n(\bm{x})=(g_1,\ldots, g_n;m_1 \ldots, m_n)$.
By the definition of $\gamma_n$, we have $m_i=m-m=0$.
Then, $\gamma_n(\bm{x}) \in D_{n}^{\rm BR_U} (X; \mathbb Z)$.

Suppose that 
\begin{align*}
\bm{x} &= (\bm{x}^{i-1} , (m_i, g_ih),  \bm{x}_{i+1})- ( \bm{x}^{i-1} , (m_i, g_i) ,  \bm{x}_{i+1})\\
&\phantom{:=} - \bigl( \bm{x}^{i-1} \uline* (m_i,g_i),( (m_i, h), \bm{x}_{i+1} ) \oline* (m_i,g_i) \bigr). 
\end{align*} 
We have
\begin{align*}
&\bigl( \bm{x}^{i-1} \uline* (m_i,g_i),\bigl( (m_i, h), \bm{x}_{i+1}\bigr) \oline* (m_i,g_i) \bigr)= \bigl((A_1, g_i^{-1}g_1g_i),\ldots,(A_{i-1}, g_i^{-1}g_{i-1}g_i),\\
&\hspace{30ex}(m_i\varphi(g_i),h),(m_{i+1}\varphi(g_i), g_{i+1}),\ldots,(m_n\varphi(g_i), g_n)
\bigr),
\end{align*}
where $A_j=m_jg_i+m_i(\varphi(g_i)-g_i)$. We note that $A_j-A_{j+1}=(m_j-m_{j+1})g_i$ for $j$ with $1\le j <i-1$ and $A_{i-1}-m_i\varphi(g_i)=(m_{i-1}-m_i)g_i$. Then,
\begin{align*}
&\gamma_n\bigl(( \bm{x}^{i-1} \uline* (m_i,g_i),\bigl( (m_i, h), \bm{x}_{i+1}\bigr) \oline* (m_i,g_i) )\bigr)\\
&=(g_i^{-1}\bm{g}^{i-1}g_i, h, \bm{g}_{i+1}; m_1'g_i, \ldots, m_{i-1}'g_i, m_i'\varphi(g_i),\ldots, m_{n-1}'\varphi(g_i),m_n\varphi(g_i)),
\end{align*}
where $m'_i:=m_i-m_{i+1}$.
Then, we have
\begin{align*}
&\gamma_n(\bm{x})=(\bm{g}^{i-1},g_i h,\bm{g}_{i+1};m_1', \ldots, m_{n-1}', m_n) -(\bm{g};m'_1, \ldots, m'_{n-1}, m_n)\\
&\phantom{\gamma_n(\bm{x})=} -(g_i^{-1}\bm{g}^{i-1}g_i, h, \bm{g}_{i+1} ; m_1'g_i, \ldots, m_{i-1}'g_i,\\
&\hspace{25ex}m_i' \varphi(g_i), \ldots, m_{n-1}' \varphi(g_i), m_n\varphi(g_i)) \in D_{n}^{\rm BR} (X; \mathbb Z).
\end{align*}
Hence $\gamma_n (D_{n}^{\rm BR} (X; \mathbb Z)) \subset D_{n}^{\rm BR_U} (X; \mathbb Z)$.

Let $\bm{z}=(\bm{g}; \bm{x}^{i-1},0, \bm{x}_{i+1} )$ be an element of $ D_{n}^{\rm BR_U} (X; \mathbb Z)$.
Put $a_k:=\sum_{j=k}^n m_j$.
Then we have $a_n=m_n$ and $a_k-a_{k+1}=m_k$ for all $1\leq k\leq n-1$. Moreover, $a_i=a_{i+1}$.
Then, we have
\[
\gamma_n\bigl((a_1,g_1), \ldots, (a_n,g_n)\bigr)=\bm{z} \in D_{n}^{\rm BR_U} (X; \mathbb Z).
\]
Let $\bm{z}= ({\bm{g}}^{i-1}, g_ih, \bm{g}_{i+1}; \bm{x} )
-(\bm{g}; \bm{x} )
-(g_i^{-1}\bm{g}^{i-1}g_i, h, \bm{g}_{i+1};\bm{x}^{i-1}g_i, \bm{x}_{i}\varphi(g_i)) \in D_{n}^{\rm BR_U} (X; \mathbb Z).$
Put $a_k:=\sum_{j=k}^n m_j$.
Then we have $a_n=m_n$ and $a_k-a_{k+1}=m_k$ for all $1\leq k\leq n-1$.
Define $x_k:=(a_k, g_k)$ and $\bm{x}:=(x_1, \ldots, x_n)$.
Then, we have
\begin{align*}
&\gamma_n\Bigl( (\bm{x}^{i-1} , (a_i, g_ih),  \bm{x}_{i+1})
- ( \bm{x}^{i-1} , (a_i, g_i) ,  \bm{x}_{i+1})\\
&\hspace{25ex}-( \bm{x}^{i-1} \uline* (a_i,g_i), ( (a_i, h), \bm{x}_{i+1} ) \oline* (a_i,g_i))\Bigr)=\bm{z}.
\end{align*}
Hence $\gamma_n^{-1}(D_{n}^{\rm BR_U} (X; \mathbb Z)) \subset D_{n}^{\rm BR} (X; \mathbb Z)$.
This completes the proof.
\end{proof}

\subsubsection{The induced homomorphism $\psi_{n,\lambda} $}

We fix a group homomorphism $\lambda:G\to A$. We also define a map $\tilde \lambda:C_n^{\rm BR_U}(X; \mathbb  Z)\to A$ as $\tilde \lambda(g_1,\ldots,g_n;\bm{m}):=\lambda(g_1)$. In addition, we define a map $\psi_{n,\lambda}:C_n^{\rm BR_U}(X; \mathbb  Z) \to C_n^{\rm gp} (M; \mathbb Z)_G\otimes A$ by
$\psi_{n,\lambda}(\bm{g};\bm{m}):=\psi_n(\bm{g};\bm{m}) \otimes \tilde{\lambda}(\bm{g};\bm{m})$ for any integer $n \ge 2$.

\begin{lemma}\label{len-brugpchain}
It holds that $\psi_{n,\lambda} (D_{n}^{\rm BR_U} (X; \mathbb Z)) \subset D_{n}^{\rm gp} (M; \mathbb Z)_G \otimes A$.
Therefore,
$\psi_{n,\lambda}$ induces the homomorphism  
$$\psi_{n,\lambda } : C_{n}^{\rm norBR_U} (X; \mathbb Z) \to C_{n}^{\rm norgp} (M; \mathbb Z)_G \otimes A.$$
\end{lemma}
\begin{proof}
It suffices to show that
$\psi_{n,\lambda}(\bm{g};\bm{m}^{i-1},0, \bm{m}_{i+1}) \in D_{n}^{\rm gp} (M; \mathbb Z)_G \otimes A$ for $i$ with $1 \le i \le n-1$ and $\psi_{n,\lambda}(\bm{z}_i) \in D_{n}^{\rm gp} (M; \mathbb Z)_G \otimes A$
for $i$ with $1 \le i \le n$, 
where
\[ \bm{z}_i:= (\bm{g}^{i-1}, g_i h, \bm{g}_{i+1};\bm{m})- (\bm{g};\bm{m} ) - (g_i^{-1}\bm{g}^{i-1}g_i, h, \bm{g}_{i+1};\bm{m}^{i-1}g_i,\bm{m}_i \varphi(g_i) ). \]
We have
\begin{align*}
&\psi_n(\bm{g};\bm{m}_{i-1}, 0, \bm{m}_{i+1}) \\
&=\displaystyle \sum_{\bm{k} \in \mathcal{K}^{n}} (-1)^{|\bm{k}|} (m_1\bm{g}_{\bm{k}_1},\ldots, m_{i-1}\bm{g}_{\bm{k}_{i-1}},0,m_{i+1} \bm{g}_{\bm{k}_{i+1}}, \ldots, m_n \bm{g}_{\bm{k}_n})\in D_{n}^{\rm gp} (M; \mathbb Z)_G.
\end{align*}
This inclusion  $\psi_n(\bm{g}; \bm{m}^{i-1},0, \bm{m}_{i+1}) \in D_{n}^{\rm gp} (M; \mathbb Z)_G$ implies $\psi_{n,\lambda}(\bm{g}; \bm{m}^{i-1},0, \bm{m}_{i+1} ) \in D_{n}^{\rm gp} (M; \mathbb Z)_G \otimes A$.

Next, we show $\psi_{n,\lambda}(\bm{z}_i) \in D_{n}^{\rm gp} (M; \mathbb Z)_G \otimes A$.
\vspace{1ex}

\noindent
For $i\not=1$, we can prove that $\psi_{n}(\bm{z}_i)=0\in D_{n}^{\rm gp} (M; \mathbb Z)_G$ as follows.
Let $s$ be an integer with $1 < s < i-1$. By the direct calculation, we have
{\normalsize
\begin{align}
&\psi_n (g_i^{-1}g_1g_i, \ldots, g_i^{-1}g_{i-1}g_i,h,g_{i+1}, \ldots, g_n; m_1g_i, \ldots, m_{i-1}g_i, m_i\varphi(g_i), \ldots, m_n\varphi(g_i)) \nonumber \\
&\begin{array}{l@{}r@{}l@{}}
= \displaystyle \sum_{{\bm k}\in \mathcal K_i^0} (-1)^{|\bm{k}|} \Bigl(  &m_1     &\varphi (g_1^{k_1}) g_2^{k_2} \cdots g_{i-1}^{k_{i-1}} g_ih^{0} g_{i+1}^{k_{i+1}}\cdots g_n^{k_n},\ldots, \\
\phantom{=} &m_{s}   &\varphi (g_1^{k_1}\cdots g_s^{k_s}) g_{s+1}^{k_{s+1}} \cdots g_{i-1}^{k_{i-1}}g_ih^{0}g_{i+1}^{k_{i+1}}\cdots g_n^{k_n},\ldots, \\
\phantom{=} &m_{i-1} &\varphi (g_1^{k_1} \cdots g_{i-1}^{k_{i-1}})g_{i}h^{0}g_{i+1}^{k_{i+1}}\cdots g_n^{k_n},\\
\phantom{=} &m_{i}   &\varphi (g_1^{k_1} \cdots g_{i-1}^{k_{i-1}}g_ih^{0})g_{i+1}^{k_{i+1}}\cdots g_n^{k_n}, \ldots, \\
\phantom{=} &m_{n}   &\varphi (g_1^{k_1} \cdots g_{i-1}^{k_{i-1}}g_{i}h^{0}g_{i+1}^{k_{i+1}}\cdots g_n^{k_n}) \Bigr)
\end{array} \label{xxx3} \\
&\begin{array}{l@{}r@{}l@{}}
\phantom{=}+\displaystyle \sum_{{\bm k}\in \mathcal K_i^1} (-1)^{|\bm{k}|} \Bigl( &m_1     &\varphi (g_1^{k_1})  g_2^{k_2} \cdots g_{i-1}^{k_{i-1}} g_ih^{1}  g_{i+1}^{k_{i+1}}\cdots g_n^{k_n}, \ldots, \\
\phantom{=} &m_{s}   &\varphi (g_1^{k_1}\cdots g_s^{k_s}) g_{s+1}^{k_{s+1}} \cdots g_{i-1}^{k_{i-1}}g_ih^{1}g_{i+1}^{k_{i+1}}\cdots g_n^{k_n}, \ldots, \\
\phantom{=} &m_{i-1} &\varphi (g_1^{k_1} \cdots g_{i-1}^{k_{i-1}})g_{i}h^{1}g_{i+1}^{k_{i+1}}\cdots g_n^{k_n}, \\
\phantom{=} &m_{i}   &\varphi (g_1^{k_1} \cdots g_{i-1}^{k_{i-1}}g_ih^{1})g_{i+1}^{k_{i+1}}\cdots g_n^{k_n}, \ldots, \\
\phantom{=} &m_{n}   &\varphi (g_1^{k_1} \cdots g_{i-1}^{k_{i-1}}g_{i}h^{1}g_{i+1}^{k_{i+1}}\cdots g_n^{k_n}) \Bigr).
\end{array} \label{xxx4} 
\end{align}
}
We also have
{\normalsize
\begin{align}
&\psi_n (g_1, \ldots, g_{i-1},g_i,g_{i+1}, \ldots, g_n ; m_1, \ldots, m_n) \nonumber \\
&\begin{array}{l@{}r@{}l@{}}
= \displaystyle \sum_{{\bm k}\in \mathcal K_i^0} (-1)^{|\bm{k}|} \Bigl(  &m_1     &\varphi (g_1^{k_1})  g_2^{k_2} \cdots g_{i-1}^{k_{i-1}} g_i^{0}g_{i+1}^{k_{i+1}}\cdots g_n^{k_n}, \ldots, \\
\phantom{=} &m_{s}   &\varphi (g_1^{k_1} \cdots g_s^{k_s}) g_{s+1}^{k_{s+1}} \cdots g_{i-1}^{k_{i-1}}g_i^{0}g_{i+1}^{k_{i+1}}\cdots g_n^{k_n}, \ldots, \\
\phantom{=} &m_{i-1} &\varphi (g_1^{k_1} \cdots g_{i-1}^{k_{i-1}} )g_i^{0}g_{i+1}^{k_{i+1}}\cdots g_n^{k_n}, \\
\phantom{=} &m_{i}   &\varphi (g_1^{k_1} \cdots g_{i-1}^{k_{i-1}}g_i^{0})g_{i+1}^{k_{i+1}}\cdots g_n^{k_n}, \ldots, \\
\phantom{=} &m_{n}   &\varphi (g_1^{k_1} \cdots g_{i-1}^{k_{i-1}}g_i^{0}g_{i+1}^{k_{i+1}}\cdots g_n^{k_n}) \Bigr)
\end{array} \label{xxx1}	\\
&\begin{array}{l@{}r@{}l@{}}
\phantom{=} + \displaystyle \sum_{{\bm k}\in \mathcal K_i^1}(-1)^{|\bm{k}|} \Bigl(&m_1     &\varphi (g_1^{k_1})  g_2^{k_2}  \cdots g_{i-1}^{k_{i-1}}g_i^{1}g_{i+1}^{k_{i+1}}\cdots g_n^{k_n}, \ldots, \\
\phantom{=} &m_{s}   &\varphi (g_1^{k_1} \cdots g_s^{k_s}) g_{s+1}^{k_{s+1}} \cdots g_{i-1}^{k_{i-1}}g_i^{1}g_{i+1}^{k_{i+1}}\cdots g_n^{k_n},\ldots, \\
\phantom{=} &m_{i-1} &\varphi (g_1^{k_1} \cdots g_{i-1}^{k_{i-1}} )g_i^{1}g_{i+1}^{k_{i+1}}\cdots g_n^{k_n},\\
\phantom{=} &m_{i}   &\varphi (g_1^{k_1} \cdots g_{i-1}^{k_{i-1}}g_i^{1})g_{i+1}^{k_{i+1}}\cdots g_n^{k_n}, \ldots, \\
\phantom{=} &m_{n}   &\varphi (g_1^{k_1} \cdots g_{i-1}^{k_{i-1}}g_i^{1}g_{i+1}^{k_{i+1}}\cdots g_n^{k_n}) \Bigr).
\end{array}\label{xxx2}
\end{align}
}

Since $\eqref{xxx2}=-\eqref{xxx3}$, we have
{\normalsize
\begin{align*}
&\psi_n (g_1, \ldots, g_{i-1}, g_i h, g_{i+1}, \ldots, g_n ; m_1, \ldots, m_n)  \\
&= \eqref{xxx4}+\eqref{xxx1}\\
&= \eqref{xxx3}+\eqref{xxx4}+\eqref{xxx1}+\eqref{xxx2}.
\end{align*}
}
Then, we have $\psi_{n}(\bm{z}_i)=0\in D_{n}^{\rm gp} (M; \mathbb Z)_G$.

Therefore, we have $\psi_{n, \lambda}(\bm{z}_i)=0\in D_{n}^{\rm gp} (M; \mathbb Z)_G \otimes A$.
\vskip 1ex

\noindent
For the case that $i=1$, we have
\begin{align*}
&\psi_n(h, g_2, \ldots, g_n;m_1\varphi(g_1), \ldots, m_n\varphi(g_1))\\
&=\displaystyle \sum_{\bm{k} \in \mathcal K_n} (-1)^{|\bm{k}|} \bigl( m_1\varphi(g_1)  \varphi(h^{k_1})g_2^{k_2}\cdots g_n^{k_n},\ldots,m_n\varphi(g_1)  \varphi(h^{k_1} g_2^{k_2} \cdots g_n^{k_n}) \bigr)\\
&=\displaystyle \sum_{\bm{k} \in \mathcal K_n} (-1)^{|\bm{k}|} \bigl( m_1  \varphi(h^{k_1})g_2^{k_2} \cdots g_n^{k_n}, \ldots, m_n  \varphi(h^{k_1} g_2^{k_2} \cdots g_n^{k_n}) \bigr)\\
&=\displaystyle \sum_{\bm{k} \in \mathcal K_n} (-1)^{|\bm{k}|} \bigl( m_1  \varphi((g_1h)^{k_1})g_2^{k_2} \cdots g_n^{k_n}, \ldots, m_n  \varphi((g_1h)^{k_1} g_2^{k_2} \cdots g_n^{k_n}) \bigr)\\
&= \psi_n(g_1h, g_2, \ldots, g_n;m_1, \ldots, m_n)\\
&= \psi_n(g_1, g_2, \ldots, g_n;m_1, \ldots, m_n)
\end{align*}
since $k_1=0$. Then, we have
\begin{align*}
&\psi_{n, \lambda}(g_1 h, g_2 ,\ldots, g_n; m_1, \ldots,  m_n )\\
&= \psi_{n}(g_1 h, g_2 ,\ldots, g_n; m_1, \ldots,  m_n )\otimes \lambda(g_1h)\\
&= \psi_{n}(g_1 h, g_2 ,\ldots, g_n; m_1, \ldots,  m_n )\otimes(\lambda(g_1)+\lambda(h))\\
&= \psi_{n}(g_1, g_2 ,\ldots, g_n; m_1, \ldots,  m_n )\otimes \lambda(g_1) \\
&\phantom{=} +\psi_n(h, g_2, \ldots, g_n ; m_1\varphi(g_1), \ldots, m_n\varphi(g_1))\otimes\lambda(h)\\
&= \psi_{n, \lambda}(g_1, g_2 ,\ldots, g_n; m_1, \ldots,  m_n ) +\psi_{n, \lambda}(h, g_2, \ldots, g_n ; m_1\varphi(g_1), \ldots, m_n\varphi(g_1)).
\end{align*}
Then, we have $\psi_{n, \lambda}(\bm{z}_i)=0 \in D_{n}^{\rm gp} (M; \mathbb Z)_G \otimes A$.
\end{proof}

\subsection{Cocycles of $G$-Alexander MCB with the trivial $X$-set}

Let $X=\sqcup_{m \in M}(\{m\} \times G)$ be the $G$-Alexander multiple conjugation biquandle of $(M, \varphi)$. 
We consider the chain complex $C_\ast (X; \mathbb Z)$ defined in Subsection~\ref{subsec_homology_MCB}, where the $X$-set $Y$ is the trivial $X$-set. 
Thus 
the first element of $C_n(X;\mathbb Z)$ is omitted. For example, an element $\langle y \rangle \langle a\rangle \langle b \rangle - \langle y \rangle \langle a,ab \rangle + \langle y \rangle \langle b,ab \rangle \in C_2(X;\mathbb Z)$ is written by $\langle a \rangle \langle b \rangle - \langle a,ab \rangle + \langle b,ab \rangle$ for simplicity. 

Define a homomorphism  ${\rm proj}_n: P_n(X; \mathbb Z) \to C_n^{\rm BR}(X; \mathbb Z)$ by
\[
{\rm proj}_n(\langle \boldsymbol{x}_1 \rangle \cdots \langle \boldsymbol{x}_k \rangle) =
\left\{
\begin{array}{ll}
(x_{11}, x_{21}, \ldots , x_{k1}) & (k=n)\\
0 & (\mbox{otherwise})
\end{array}
\right.
\]
for $n \ge 2$,
where $\langle \boldsymbol{x}_j \rangle $ means $\langle {x}_{j1} , x_{j2}, \ldots , {x}_{jn_j}  \rangle $  for each $j\in\{ 1, \ldots , k\}$.
We define ${\rm proj}_n=0$ for  $n \le 1$.

\begin{lemma}
It holds that ${\rm proj}_n(D_n(X; \mathbb Z)) \subset D_n^{\rm BR}(X; \mathbb Z)$.
\end{lemma}

\begin{proof}
Let $(m_1, g_1), \ldots , (m_n, g_n)$ be elements of $X$ such that $m_i=m_{i+1}$ for some $i$.
Put $x_i:=(m_i, g_j) $. 
We show that
\begin{align*}
&{\rm proj}_n \Bigl( \langle x_1 \rangle \cdots \langle x_{i-1} \rangle \langle x_i \rangle \langle x_{i+1} \rangle \langle x_{i+2} \rangle \cdots \langle x_n \rangle - \\
&\hspace{6.6ex}\langle x_1 \rangle \cdots \langle x_{i-1} \rangle \langle \langle x_i \rangle \langle x_{i+1} \rangle \rangle  \langle x_{i+2} \rangle\cdots \langle x_n \rangle \Bigr) \in D_n^{\rm BR}(X; \mathbb Z).
\end{align*}
We have
\begin{align*}
&{\rm proj}_n(
\langle x_1 \rangle \cdots \langle x_{i-1}\rangle \langle x_i \rangle \langle x_{i+1} \rangle \langle x_{i+2}\rangle \cdots \langle x_n \rangle\\
&\phantom{=} -\langle x_1 \rangle \cdots \langle x_{i-1}\rangle \langle \langle x_i \rangle \langle x_{i+1}\rangle \rangle \langle x_{i+2}\rangle \cdots \langle x_n \rangle)\\
&={\rm proj}_n( \langle x_1 \rangle \cdots \langle x_{i-1}\rangle \langle x_i \rangle \langle x_{i+1} \rangle \langle x_{i+2}\rangle \cdots \langle x_n \rangle\\
&\phantom{=} -\langle x_1 \rangle \cdots \langle x_{i-1}\rangle (\langle  x_i, x_i x_{i+1} \rangle - \langle x_{i+1}, x_i x_{i+1} \rangle ) \langle x_{i+2}\rangle \cdots \langle x_n \rangle )\\
&={\rm proj}_n (\langle x_1 \rangle \cdots \langle x_{i-1} \rangle \langle x_i \rangle \langle x_{i+1}\rangle \langle x_{i+2} \rangle \cdots \langle x_n \rangle) \\
&\phantom{=} -{\rm proj}_n( \langle x_1 \rangle \cdots \langle x_{i-1} \rangle \langle x_i, x_i x_{i+1} \rangle \langle x_{i+2} \rangle \cdots \langle x_n\rangle) \\
&\phantom{=} +{\rm proj}_n(\langle x_1 \rangle \cdots \langle x_{i-1} \rangle \langle x_{i+1},x_i x_{i+1}\rangle \langle x_{i+2} \rangle \cdots \langle x_n \rangle )\\
&= {\rm proj}_n ( \langle x_1 \rangle \cdots \langle x_{i-1} \rangle \langle x_i \rangle  \langle x_{i+1} \rangle \langle x_{i+2} \rangle \cdots \langle x_n\rangle )\\
&={\rm proj}_n ( \langle(m_1,g_1) \rangle \cdots \langle (m_i,g_i) \rangle \langle(m_i,g_{i+1}) \rangle \cdots \langle (m_n,g_n) \rangle )\\
&=((m_1,g_1),\ldots,(m_i,g_i),(m_i,g_{i+1}),\ldots,(m_n,g_n)) \in D_n^{\rm BR}(X; \mathbb Z).
\end{align*}
We have ${\rm proj}_n(D_n(X; \mathbb Z )) \subset D_n^{\rm BR}(X; \mathbb Z)$.
\end{proof}
Therefore, the homomorphism ${\rm proj}_n: P_n(X; \mathbb Z ) \to C_n^{\rm BR}(X; \mathbb Z )$ induces the homomorphism ${\rm proj}_n: C_n(X; \mathbb Z ) \to C_n^{\rm norBR}(X; \mathbb Z )$.
\begin{lemma}\label{lemma:projA}
The map ${\rm proj}: C_*(X; \mathbb Z ) \to C_*^{\rm norBR}(X; \mathbb Z )$ is a chain map, that is, it holds that for any integer $n$
\[
{\rm proj}_{n-1} \circ \partial_n =\partial_n^{\rm norBR} \circ {\rm proj}_n.
\]
\end{lemma}
\begin{proof}
It is sufficient to consider $\langle x_1 \rangle \cdots \langle x_n \rangle $ 
in $C_n(X; \mathbb Z)$. For $\langle x_1 \rangle \cdots \langle x_n\rangle $, we have
\begin{align*}
&\partial_n^{\rm norBR} \circ {\rm proj}_n(\langle x_1 \rangle \cdots \langle x_n \rangle)\\
&=\partial_n^{\rm norBR}(x_1,\ldots,x_n)\\
&=\displaystyle\sum^{n}_{i=1}(-1)^{i-1}\{ (\bm{x}^{i-1}, \bm{x}_{i+1})-(\bm{x}^{i-1} \uline{*}x_i, \bm{x}_{i+1} \oline{*}x_i) \}\\
&=\displaystyle\sum^{n}_{i=1} (-1)^{i-1} 
\bigl\{ 
{\rm proj}_{n-1}( \langle x_1 \rangle \cdots \langle x_{i-1}\rangle \langle x_{i+1}\rangle \cdots \langle x_n \rangle )\\
&\phantom{=}
-{\rm proj}_{n-1}
( \langle  x_1 \uline{*} x_i  \rangle \cdots\langle x_{i-1} \uline{*} x_i  \rangle  \langle x_{i+1}  \oline{*}  x_i  \rangle \cdots \langle x_n  \oline{*}  x_i  \rangle ) 
\bigr\} \\
&= {\rm proj}_{n-1} \circ \partial_n (\langle x_1 \rangle \cdots \langle x_n \rangle ).
\end{align*}
We have $ {\rm proj}_{n-1} \circ \partial_n =\partial_n^{\rm norBR} \circ {\rm proj}_n$.
\end{proof}

Fix a group homomorphism $\lambda: G \to A$. 
As a consequence of Lemma~\ref{lemma:projA} and Subsection~\ref{subsection:A3}, we have the sequence
\[
C_n (X; \mathbb Z) \overset{{\rm proj}_n}{\longrightarrow} C_n^{\rm norBR} (X; \mathbb Z) \overset{\gamma_n}{\longrightarrow} C_{n}^{\rm norBR_U} (X; \mathbb Z) \overset{\psi_{n,\lambda}}{ \longrightarrow } C_{n}^{\rm norgp} (M; \mathbb Z)_G\otimes A
\]
of chain groups. Therefore, we have the following theorem. 

\begin{theorem}\label{mcb_cocycles_no_X-set_main}
For any  $n$-cocycle $f: C_{n}^{\rm norgp} (M; \mathbb Z)_G \to A$, the map
$$\Phi_{f,\lambda}:=(f \otimes {\rm id}_A) \circ  \psi_{n,\lambda} \circ \gamma_n \circ  {\rm proj}_n : C_n (X; \mathbb Z)  \to A$$  
is an $n$-cocycle of the $G$-Alexander MCB $X=\sqcup_{m \in M}(\{x\} \times G)$.
\end{theorem}

The following theorem follows from the direct calculation.

\begin{theorem}\label{mcb_cocycle_no_X-setexplicit}
\begin{itemize}
\item[(1)] 
Let $f: M^2 \to A$ be a $G$-invariant $A$-multilinear map. 
The $2$-cocycle $\Phi_{f,\lambda}= (f \otimes {\rm id}_A) \circ  \psi_{2, \lambda} \circ \gamma_2  \circ {\rm proj}_2 : C_2 (X; \mathbb Z)  \to A$ of the $G$-Alexander MCB $X=\bigsqcup_{m \in M} (\{m\} \times G)$ is formulated as 
\[
\Phi_{f,\lambda} \left( \langle (m_1,g_1) \rangle \langle (m_2,g_2) \rangle \right)
=f\Bigl(m_1-m_2, m_2(1-\varphi(g_2)g_2^{-1})\Bigr) \otimes \lambda (g_1)
\]
for $\langle (m_1,g_1) \rangle \langle (m_2,g_2) \rangle \in X^2 \subset C_2 (X; \mathbb Z)$.
\item[(2)] 
Let $f: M^3 \to A$ be a $G$-invariant $A$-multilinear map. 
The $3$-cocycle $\Phi_{f,\lambda} = (f \otimes {\rm id}_A) \circ  \psi_{3, \lambda} \circ \gamma_3 \circ {\rm proj}_3 : C_3 (X; \mathbb Z)  \to A$ of the $G$-Alexander MCB $X=\bigsqcup_{m \in M} (\{m\} \times G)$ is
 formulated as 
\begin{align*}
&\Phi_{f,\lambda}( \langle (m_1,g_1) \rangle  \langle (m_2,g_2) \rangle \langle (m_3,g_3) \rangle )\\
&=f\Bigl((m_1-m_2)(1-\varphi(g_2)^{-1}g_2),m_2-m_3,m_3(1-\varphi(g_3)g_3^{-1})\Bigr) \otimes \lambda (g_1)
\end{align*}
for $\langle (m_1,g_1) \rangle  \langle (m_2,g_2) \rangle \langle (m_3,g_3) \rangle  \in X^3\subset C_3 (X; \mathbb Z)$.
\end{itemize}
\end{theorem}

\section{Cocycles of $G$-Alexander multiple conjugation biquandles with the $X$-set $X$}\label{app:22}

Throughout this section, let $X=\bigsqcup_{m \in M}(\{m\} \times G)=M \times G$
the $G$-Alexander multiple conjugation biquandle of $(M, \varphi)$.
We assume that the $X$-set $Y$ is $X$ itself. 
Our goal in this section is to give Theorem~\ref{mcb_cocycle_with_X-set_main}.

\subsection{Degenerate subcomplexes
$D^{\rm BR}_*(X;\mathbb Z)_X$, $D^{\rm BR_U}_*(X;\mathbb Z)_X$ and the induced homomorphisms $\gamma_n$, $\psi_n$ and $\psi_{n,\lambda}$}\label{subsection:B3}

\subsubsection{The degenerate subcomplex $D^{\rm BR}_* (X;\mathbb Z)_X$ of $C^{\rm BR}_* (X;\mathbb Z)_X$}

Let $D^{\rm BR}_n(X;\mathbb Z)_X$ be the subgroup of $C^{\rm BR}_n(X;\mathbb Z)_X$ generated by the elements of the following sets
\begin{align*}
&
\displaystyle \bigcup_{ i=1}^{ n-1} \Big\{ ( \bm{x}^{i-1} , (m, g),  (m, h),  \bm{x}_{i+2}) ~\Big|~\bm{x} \in X \times X^n,~m \in M,~g,h \in G \Big\} \text{ and}\\
&
\displaystyle \bigcup_{i=1}^{n} \left\{
\begin{array}{lcl}
\begin{minipage}{6.6cm}
{\normalsize
$\hspace{1.8ex}(\bm{x}^{i-1} , (m, gh),  \bm{x}_{i+1} ) - ( \bm{x}^{i-1} , (m, g),  \bm{x}_{i+1})$\\
$ - ( \bm{x}^{i-1} \underline{*} (m,g),~ \bigl( (m, h),  \bm{x}_{i+1} \bigr) \overline{*} (m,g) )$}
\end{minipage}
& \Bigg| &
\begin{minipage}{2.7cm}
{\normalsize
$\bm{x} \in X \times X^n,\\ m \in M,~g,h \in G$
}
\end{minipage}
\end{array}
\right\}
\end{align*}
for   $n \ge 2$,
where we write
\begin{align*}
&\left( \bm{x}^{i-1} \underline{*} (m,g),~\bigl( (m, h),  \bm{x}_{i+1} \bigr) \overline{*} (m,g) \right)\\
&=(x_0 \uline*(m,g), \ldots, x_{i-1}\uline*(m,g), (m,h)\oline*(m,g),~x_{i+1}\oline*(m,g), \ldots, x_n\oline*(m,g)).
\end{align*}
Define $D^{\rm BR}_n(X;\mathbb Z)_X:=0$ for  $n \le 1$.

\begin{lemma} 
$D^{\rm BR}_*(X;\mathbb Z)_X:=(D^{\rm BR}_n(X;\mathbb Z)_X, \partial_n^{\rm BR})_{n \in \mathbb Z}$ is a subcomplex of $C^{\rm BR}_*(X;\mathbb Z)_X$.
\end{lemma}

The chain complex 
$$C^{\rm norBR}_* (X;\mathbb Z)_X:=C^{\rm BR}_* (X;\mathbb Z)_X/D^{\rm BR}_* (X;\mathbb Z)_X$$ 
 determines the  homology group $H^{\rm nor BR}_n(X;\mathbb Z)_X$. 
In the ordinary way, for an abelian group $A$, we have the (co)homology theory with the coefficient group $A$ and the homology group $H^{\rm nor BR}_n(X;A)_X$ and the cohomology group $H_{\rm nor BR}^n(X;A)_X$ are defined.

\subsubsection{The degenerate subcomplex $D^{\rm BR_U}_*(X;\mathbb Z)_X$ of $C^{\rm BR_U}_* (X;\mathbb Z)_X$}

Let $D^{\rm BR_U}_n(X;\mathbb Z)_X$ be the subgroup of $C^{\rm BR_U}_n(X;\mathbb Z)_X$ generated by the elements of the following sets
\begin{align*}
\displaystyle \bigcup_{ i=1}^{ n-1} &\Big\{ (\bm{g};\bm{m}^{i-1}, 0, \bm{m}_{i+1}) ~\Big|~ \bm{g} \in G \times G^n, \bm{m} \in M \times M^n\Big\} \mbox{ and}\\
\displaystyle \bigcup_{i=1}^{n} &\Bigl\{
\begin{array}{lcl}
\begin{minipage}{6.2cm}
\hspace{2ex}$(\bm{g}^{i-1}, g_i h, \bm{g}_{i+1} ;\bm{m} )- (\bm{g};\bm{m} )$\\
$-(g_i^{-1} \bm{g}^{i-1}g_i, h, \bm{g}_{i+1}; \bm{m}^{i-1}g_i,~\bm{m}_i \varphi(g_i))$
\end{minipage}
& \Big| &
\begin{minipage}{3.4cm}
$h \in G$,
$\bm{g} \in G \times G^n$,\\
$\bm{m} \in M \times M^n$
\end{minipage}
\end{array}
\Bigr\}
\end{align*}
for $n \ge 2$, where we write
\begin{align*}
&(\bm{g};\bm{m}^{i-1},0,\bm{m}_{i+1}):=(g_0, g_1, \ldots, g_n;m_0, m_1, \ldots, m_{i-1}, 0, m_{i+1}, \ldots, m_n),\\
&(\bm{g}^{i-1},g_ih,\bm{g}_{i+1};\bm{m}):=(g_0, g_1, \ldots, g_{i-1},g_ih, g_{i+1},\ldots ,g_n;m_0, m_1, \ldots, m_n) \text{ and}\\
&(g_i^{-1}\bm{g}^{i-1}g_i, h, \bm{g}_{i+1} ;\bm{m}^{i-1}g_i, \bm{m}_i \varphi(g_i) ):=\bigl( (g_{i}^{-1}g_0g_i), (g_{i}^{-1}g_1g_i), \ldots, (g_i^{-1}g_{i-1}g_i),\\
&\hspace{11ex}h,g_{i+1},\ldots,g_n;m_0g_i, m_1g_i, \ldots, m_{i-1}g_i,m_i\varphi(g_i), \ldots, m_n\varphi(g_i) \bigr).
\end{align*}
Define $D^{\rm BR_U}_n(X;\mathbb Z)_X:=0$ for  $n \le 1$.

\begin{lemma}
$D^{\rm BR_U}_*(X;\mathbb Z)_X:=(D^{\rm BR_U}_n (X;\mathbb Z)_X, \partial_n^{\rm BR_U})_{n \in \mathbb Z}$ is a subcomplex of $C^{\rm BR_U}_*(X;\mathbb Z)_X$.
\end{lemma}

The chain complex 
$$C^{\rm norBR_U}_*(X;\mathbb Z)_X:=C^{\rm BR_U}_*(X;\mathbb Z)_X/D^{\rm BR_U}_* (X;\mathbb Z)_X$$ 
determines the homology group $H^{\rm nor BR_U}_n(X;\mathbb Z)_X$.  
In the ordinary way, for an abelian group $A$, we have the (co)homology theory with the coefficient group $A$ and the homology group $H^{\rm nor BR_U}_n(X;A)_X$ and the cohomology group $H_{\rm nor BR_U}^n(X;A)_X$ are defined.  

\subsubsection{The induced homomorphism $\gamma_n$}

Next lemma shows that the isomorphism  
$\gamma_n: C_n^{\rm BR} (X; \mathbb Z)_X \to C_{n}^{\rm BR_U} (X; \mathbb Z)_X$
defined in Subsection~\ref{subsection:B1} induces the isomorphism  $\gamma_n:C_n^{\rm norBR} (X; \mathbb Z)_X \to C_{n}^{\rm norBR_U} (X; \mathbb Z)_X$, where we denote it by the same symbol  $\gamma_n$ for simplicity.

\begin{lemma}\label{4lem-resgam} 
It holds that $\gamma_n(  D_n^{\rm BR} (X; \mathbb Z)_X ) = D_{n}^{\rm BR_U} (X; \mathbb Z)_X$.  
Therefore 
$\gamma_n$ induces the isomorphism
$$\gamma_n: C_n^{\rm norBR} (X; \mathbb Z)_X \to C_{n}^{\rm norBR_U} (X; \mathbb Z)_X.$$
\end{lemma}
\subsubsection{The induced homomorphisms $\psi_{n}$ and $\psi_{n,\lambda}$}

Next lemma shows that the map  $\psi_{n} : C_{n}^{\rm BR_U} (X; \mathbb Z)_X \to C_{n+1}^{\rm gp} (M; \mathbb Z)_G$
defined in Subsection~\ref{subsection:B1} induces the homomorphism $\psi_{n} : C_{n}^{\rm norBR_U} (X; \mathbb Z)_X \to C_{n+1}^{\rm norgp} (M; \mathbb Z)_G$, where we denote it by the same symbol  $\psi_n$ for simplicity.
\begin{lemma}\label{lem-resgam}
It holds that
$\psi_n(D_{n}^{\rm BR_U} (X; \mathbb Z)_X)\subset D_{n+1}^{\rm gp} (M; \mathbb Z)_G$.
Therefore 
$\psi_n$ induces the homomorphism
\[ \psi_{n} : C_{n}^{\rm norBR_U} (X; \mathbb Z)_X \to C_{n+1}^{\rm norgp} (M; \mathbb Z)_G. \]
\end{lemma}
In addition, a homomorphism $\psi_{n,\lambda} : C_{n}^{\rm norBR_U} (X; \mathbb Z)_X \to C_{n+1}^{\rm norgp} (M; \mathbb Z)_G \otimes A$ is also defined as follows:
We fix a group homomorphism $\lambda:G\to A$. 
We define a map $\tilde \lambda:C_n^{\rm BR_U}(X; \mathbb  Z)_X\to A$ as $\tilde \lambda(g_0, g_1,\ldots,g_n; m_0 , m_1,\ldots,m_n)=\lambda(g_0)$. 
Define a map $\psi_{n,\lambda}:C_n^{\rm BR_U}(X; \mathbb  Z)_X \to C_{n+1}^{\rm gp} (M; \mathbb Z)_G\otimes A$ by 
$\psi_{n,\lambda}= \psi_n \otimes \tilde{\lambda} $ for $n \geq 1$. Define $\psi_{n,\lambda}:=0$ for $n <1$. 

\begin{lemma}\label{len-brugpchain}
It holds that $\psi_{n,\lambda} (D_{n}^{\rm BR_U} (X; \mathbb Z)_X) \subset D_{n+1}^{\rm gp} (M; \mathbb Z)_G \otimes A$.
Therefore 
$\psi_{n,\lambda}$ induces the homomorphism  
\[ \psi_{n,\lambda } : C_{n}^{\rm norBR_U} (X; \mathbb Z)_X \to C_{n+1}^{\rm norgp} (M; \mathbb Z)_G \otimes A.\]
\end{lemma}
\subsection{Cocycles of $G$-Alexander MCB with the $X$-set $X$}
Define a homomorphism  ${\rm proj}_n: P_n(X;\mathbb Z)_X \to C_n^{\rm BR}(X; \mathbb Z )_X$ by
\[ {\rm proj}_n(\langle x_0 \rangle \langle \bm{x}_1 \rangle \cdots \langle \bm{x}_k \rangle) = \left\{
\begin{array}{ll}
(x_0, x_{11},\ldots , x_{k1}) & (k=n)\\
0 & (\mbox{otherwise})
\end{array}
\right.
\]
for  $n \ge 2$,
where $\langle \boldsymbol{x}_j \rangle $ means $\langle {x}_{j1} , x_{j2}, \ldots , {x}_{jn_j}  \rangle $  for each $j\in\{ 1, \ldots , k\}$.
Define ${\rm proj}_n=0$ for  $n \le 1$.

\begin{lemma}
It holds that ${\rm proj}_n(D_n(X; \mathbb Z )_X) \subset D_n^{\rm BR}(X; \mathbb Z )_X$.
\end{lemma}

Therefore, the homomorphism ${\rm proj}_n: P_n(X; \mathbb Z )_X \to C_n^{\rm BR}(X;\mathbb Z )_X$ induces the homomorphism ${\rm proj}_n: C_n(X; \mathbb Z )_X \to C_n^{\rm norBR}(X; \mathbb Z )_X$.
\begin{lemma}\label{lemma:projB}
The map ${\rm proj}$ is a chain map, that is, it holds that 
\[
{\rm proj}_{n-1} \circ \partial_n =\partial_n^{\rm norBR} \circ {\rm proj}_n.
\]
\end{lemma}
As a consequence of Lemma \ref{lemma:projB} and Subsection \ref{subsection:B3}, we have two sequences
\begin{align*}
&
C_n (X; \mathbb Z)_X \overset{{\rm proj}_n}{\longrightarrow} C_n^{\rm norBR} (X; \mathbb Z)_X \overset{\gamma_n}{\longrightarrow} C_{n}^{\rm norBR_U} (X; \mathbb Z)_X \overset{\psi_{n}}{ \longrightarrow } C_{n+1}^{\rm norgp} (M; \mathbb Z)_G \otimes A,\\
&
C_n (X; \mathbb Z)_X \overset{{\rm proj}_n}{\longrightarrow} C_n^{\rm norBR} (X; \mathbb Z)_X \overset{\gamma_n}{\longrightarrow} C_{n}^{\rm norBR_U} (X; \mathbb Z)_X \overset{\psi_{n,\lambda}}{ \longrightarrow } C_{n+1}^{\rm norgp} (M; \mathbb Z)_G\otimes A.
\end{align*}
Therefore, we have the following theorem.
\begin{theorem}\label{mcb_cocycle_with_X-set_main}  
For any $(n+1)$-cocycle $f: C_{n+1}^{\rm norgp} (M; \mathbb Z)_G \to A$, the maps
\begin{align*}
\Phi_{f} &:= f  \circ  \psi_{n} \circ \gamma_n \circ  {\rm proj}_n : C_n (X; \mathbb Z)_X  \to A \text{ and}\\
\Phi_{f,\lambda} &:=(f \otimes {\rm id}_A) \circ  \psi_{n,\lambda} \circ \gamma_n \circ  {\rm proj}_n : C_n (X; \mathbb Z)_X  \to A
\end{align*}
are $n$-cocycles  of the $G$-Alexander MCB $X=\bigsqcup_{m \in M} (\{m\} \times G)$.
\end{theorem}

\begin{remark}
The cocycle $\Phi_{f}$ does not use the information $g_0$ of $(m_0, g_0)\in X=M\times G$. Therefore, we may replace the $X$-set $X=M\times G$ with $M$, that is, $\Phi_{f}$ can be also defined as the map from $C_n (X; \mathbb Z)_X $ to $A$.
\end{remark}

\begin{theorem}\label{mcb_cocycle_with_X-set_main_explicit1}
\begin{itemize}
\item[(1)] 
Let $f: M^3 \to A$ be a $G$-invariant $A$-multilinear map. 
The $2$-cocycle $\Phi_f= f \circ  \psi_2 \circ \gamma_2   : C_2 (X; \mathbb Z)_X  \to A$ of the $G$-Alexander MCB $X=M \times G=\bigsqcup_{m \in M} (\{m\} \times G)$ is formulated as 
\begin{align*}
&
\hspace{-2ex}\Phi_f \bigl( \langle (m_0 , g_0 ) \rangle \langle (m_1,g_1) \rangle  \langle (m_2,g_2) \rangle \bigr)\\
&
= f \bigl( m_0'(1-\varphi(g_1)^{-1}g_1) ,m_1',m_2(1-\varphi(g_2)g_2^{-1}) \bigr)
\end{align*}
for any $\langle (m_0 , g_0 ) \rangle \langle (m_1,g_1) \rangle  \langle (m_2,g_2) \rangle  \in X\times X^2$,
where $m_i':=m_i -m_{i+1}$.
\item[(2)] 
Let $f: M^4 \to A$ be a $G$-invariant $A$-multilinear map. 
The $3$-cocycle $\Phi_f= f \circ  \psi_3 \circ \gamma_3   : C_3 (X; \mathbb Z)_X  \to A$ of the $G$-Alexander MCB $X=M \times G=\bigsqcup_{m \in M} (\{m\} \times G)$ is
 formulated as 
\begin{align*}
&
\hspace{-2ex}\Phi_f \bigl( \langle (m_0,g_0) \rangle \langle ( m_1,g_1 ) \rangle \langle (m_2,g_2) \rangle  \langle (m_3,g_3) \rangle \bigr) \\
&
=f\bigl( m_0' (1-\varphi(g_1)^{-1}g_1),m_1', m_2', m_3(1-\varphi(g_3)g_3^{-1}) \bigr)\\
&
\phantom{=} -f\bigl( m_0' (1-\varphi(g_1)^{-1}g_1)g_2, m_1' g_2, m_2' \varphi(g_2), m_3(1-\varphi(g_3)g_3^{-1})\varphi(g_2) \bigr)
\end{align*}
for any $ \langle (m_0,g_0) \rangle  \langle ( m_1,g_1 ) \rangle \langle (m_2,g_2) \rangle  \langle (m_3,g_3) \rangle \in X \times X^3$,
where $m_i':=m_i -m_{i+1}$.
\end{itemize}
\end{theorem}

\begin{theorem}\label{mcb_cocycle_with_X-set_main_explicit2}
\begin{itemize}
\item[(1)] 
Let $f: M^3 \to A$ be a $G$-invariant $A$-multilinear map. 
The $2$-cocycle $\Phi_{f, \lambda}= (f \otimes {\rm id}_A) \circ  \psi_{2, \lambda} \circ \gamma_2   : C_2 (X; \mathbb Z)_X  \to A$ of the $G$-Alexander MCB $X=\sqcup_{m \in M} (\{x\} \times G)$ is formulated as 
\begin{align*}
&
\Phi_{f, \lambda}(\langle (m_0, g_0) \rangle \langle (m_1,g_1) \rangle  \langle (m_2,g_2) \rangle )\\
&
=f \Bigl( m_0' (1-\varphi(g_1)^{-1}g_1),m_1',m_2(1-\varphi(g_2)g_2^{-1}) \Bigr) \otimes \lambda (g_0)
\end{align*}
for any $\langle (m_0, g_0) \rangle \langle (m_1,g_1) \rangle  \langle (m_2,g_2) \rangle  \in X\times X^2$,
where $m_i':=m_i -m_{i+1}$.
\item[(2)] 
Let $f: M^4 \to A$ be a $G$-invariant $A$-multilinear map. 
The $3$-cocycle $\Phi_{f, \lambda}= (f \otimes {\rm id}_A) \circ  \psi_{3, \lambda} \circ \gamma_3   : C_3(X; \mathbb Z)_X  \to A$ of the $G$-Alexander MCB $X=\sqcup_{m \in M} (\{x\} \times G)$ is
 formulated as 
\begin{align*}
&
\hspace{-2.5ex}
\Phi_{f, \lambda}( \langle (m_0, g_0) \rangle   \langle ( m_1,g_1 ) \rangle  \langle (m_2,g_2) \rangle  \langle (m_3,g_3) \rangle )\\
&
\hspace{-2ex} =f \Bigl( m_0' (1-\varphi(g_1)^{-1}g_1),m_1' ,m_2', m_3(1-\varphi(g_3)g_3^{-1}) \Bigr) \otimes \lambda (g_0)\\
&
\hspace{-2ex} \phantom{=} -f \Bigl( m_0' (1-\varphi(g_1)^{-1}g_1)g_2, m_1' g_2, m_2' \varphi(g_2), m_3(1-\varphi(g_3)g_3^{-1})\varphi(g_2) \Bigr) \otimes  \lambda (g_0)
\end{align*}
for any $ \langle (m_0, g_0) \rangle  \langle ( m_1,g_1 ) \rangle \langle (m_2,g_2) \rangle   \langle (m_3,g_3) \rangle  \in X \times X^3$,
where $m_i':=m_i -m_{i+1}$.
\end{itemize}
\end{theorem}

\end{document}